\documentclass{amsart}
\usepackage{amsmath, amsthm, amsfonts, amssymb}
\usepackage{hyperref}
\hypersetup{
    colorlinks=true,
    linkcolor=black,
    citecolor=black,
    urlcolor=black,
    pdfborder={0 0 0}  
}
\usepackage{mathrsfs,comment}
\usepackage{color}
\usepackage{BOONDOX-cal}
\usepackage{cases}
\usepackage{caption}
\usepackage{subcaption} 
\usepackage{pifont}
\usepackage{tikz}
\usetikzlibrary{patterns} 
\usepackage{tikzit}
\usepackage{array}
\usepackage{multirow}
\usepackage{longtable}
\usepackage{makecell}
\numberwithin{equation}{section} 
\newtheorem{thm}{Theorem}[section]
\newtheorem{lem}{Lemma}[section]
\newtheorem{prop}{Proposition}[section]
\newtheorem{rem}{Remark}[section]
\newtheorem{cor}{Corollary}[section]
\newtheorem{dfn}{Definition}[section]

\newtheorem{con}{Conjecture}[section]

\allowdisplaybreaks[4]

\tikzstyle{Blue Dot}=[fill=blue, draw=blue, shape=circle]
\tikzstyle{Yellow Dot}=[fill=yellow, draw=yellow, shape=circle]
\tikzstyle{Black}=[fill=black, draw=black, shape=circle]
\tikzstyle{Red}=[fill=red, draw=red, shape=circle]
\tikzstyle{new style 0}=[fill=white, draw=black, shape=circle]

\tikzstyle{new edge style 0}=[-, dashed, dash pattern=on 2mm off 2mm]
\tikzstyle{new edge style 1}=[-, draw=red]
\tikzstyle{new edge style 2}=[-, fill=none, draw=blue]
\tikzstyle{new edge style 3}=[<->, dashed, dash pattern=on 2mm off 2mm]

\begin{document}
\title[Fractional Triangular Maximal Operators on Locally Finite Trees]{Fractional Triangular Maximal Operators on Locally Finite Trees}
\author{Chenghao Xiang and Fayou Zhao$^\ast$ }
\address[Chenghao Xiang]{Department of Mathematics, Shanghai University, Shanghai
200444, P.R. China}
\email{chenghaoxiang@shu.edu.cn}
\address[Fayou Zhao]{Newtouch Center for Mathematics of Shanghai University,  Department of Mathematics, Shanghai University,
Shanghai 200444, P.R. China}
\email{fyzhao@shu.edu.cn}
\thanks{The research was supported by National Natural
Science Foundation of China (Grant No. 11971295).\\
$^\ast$ Corresponding author}

\begin{abstract}
We introduce and study the centered and uncentered fractional triangular maximal operators, $\mathcal{N}^\gamma$ and $\widetilde{\mathcal{N}}^\gamma$, on locally finite trees. A key geometric tool that consists of a family of overlapping triangles with total measure comparable to a ball is established, which facilitates the analysis of ball-type maximal operators. Our main results include: complete pointwise comparisons among these and classical fractional maximal operators on homogeneous trees; sufficient geometric conditions and counterexamples for general trees; and sharp $\ell^p \to \ell^q$ boundedness theorems for the new maximal operators, with optimal ranges of exponents verified by counterexamples. Finally, we compare the unboundedness regions of the ball-type maximal operator $\widetilde{\mathcal{M}}^\gamma$ and the triangular maximal operator $\widetilde{\mathcal{N}}^\gamma$ on a prototypical tree, proving that the former’s region strictly contains the latter’s for all $0 < \gamma \leq 1$.
\end{abstract}

\maketitle
\section{Introduction}\label{s1}
Let $(X, d, \mu)$ be a metric measure space and let $f \colon X \to \mathbb{R}$ be a locally integrable function. 
The centered and uncentered fractional Hardy--Littlewood maximal operators of order $\gamma > 0$ are defined by
\[
M^{\gamma}f(x) := \sup_{r>0} \frac{1}{\mu\!\bigl(\mathcal{B}(x,r)\bigr)^{\gamma}} \int_{\mathcal{B}(x,r)} \lvert f(y) \rvert \, d\mu(y), \quad x\in X,
\]
and
\[
\widetilde{M}^{\gamma}f(x) := \sup_{\mathcal{B} \ni x } \frac{1}{\mu(\mathcal{B})^{\gamma}} \int_{\mathcal{B}} \lvert f(y) \rvert \, d\mu(y), \quad x\in X.
\]
Here $\mathcal{B}(x,r)$ denotes the open ball centred at $x$ with radius $r$, and the second supremum is taken over all balls $\mathcal{B}$ in $X$ that contain $x$.  
When $\gamma = 1$, these reduce to the classical centered and uncentered Hardy--Littlewood maximal operators $M$ and $\widetilde{M}$, respectively.

The measure $\mu$ is said to be doubling if there exists a constant $C_D > 0$ such that for all $x \in X$ and $r > 0$,
\[
\mu(\mathcal{B}(x, 2r)) \leq C_D \, \mu(\mathcal{B}(x, r)).
\]
A metric measure space $(X, d, \mu)$ endowed with a doubling measure $\mu$ is called a space of homogeneous type in the sense of Coifman and Weiss \cite{CW}. 
Using a Vitali-type covering lemma adapted to the homogeneous setting, Coifman and Weiss \cite{CW} proved that the Hardy--Littlewood maximal operator $M$ is of weak type $(1,1)$ and bounded on $L^p(X)$ for every $p \in (1, \infty]$.
In contrast, for measures that do not satisfy the doubling condition, the boundedness properties of the uncentered maximal operator become considerably more delicate. Let $\mu$ be a non-negative Borel measure on $\mathbb{R}^n$, $n \geq 1$.
It is well known that the centered maximal operator $M$ is of weak type $(1,1)$ and strong type $(p,p)$ for $p>1$, as can be shown via the Besicovitch covering lemma.
However, the uncentered maximal operator $\widetilde{M}$ does not, in general, share these mapping properties when $n \geq 2$.
Indeed, as shown in \cite{SJ}, for $n=2$ there exists a Gaussian measure $\mu$ such that $\widetilde{M}$ fails to be of weak type $(1,1)$ while in one dimension $\widetilde{M}$ retains the same boundedness as $M$. These facts are closely related to the doubling condition: if $\mu$ is doubling, then $\widetilde{M}$ is of weak type $(1,1)$; otherwise, the boundedness holds only in the one-dimensional setting. Another example is that Str\"omberg \cite{STR} proved that $M$ is bounded on $L^p$ for all $p>1$ and is of weak type $(1,1)$ on symmetric spaces of the noncompact type while Ionescu \cite{I} proved that $\widetilde{M}$ is bounded on $L^p$ if and only if $p>2$.

In the present work, we turn to the setting of graphs, specifically trees, where the natural counting measure is typically non-doubling. This provides a natural framework to further investigate the behavior of Hardy--Littlewood maximal operators beyond the doubling regime. We focus on the fractional Hardy--Littlewood maximal operator on trees. To this end, we first introduce some basic concepts and notations related to trees. 

Let $\mathfrak{T}=(V_\mathfrak{T},E_\mathfrak{T})$ be a locally finite tree, i.e., a simple, connected, undirected, infinite graph with no loops in which every vertex $x$ has a finite number $\nu(x)\geq3$ of neighbors;
$\nu(x)$ is called the degree of $x$.
$V_\mathfrak{T}$ denotes the set of vertices and $E_\mathfrak{T}$ denotes the set of edges between them.
Here, a graph is simple if it has no parallel edges and there are no edges connecting one vertex and itself. 
Two distinct vertices $x$, $y$ in $\mathfrak{T}$ are neighbors if and only if they are connected by exactly one edge, denoted by $x\sim y$.
$\mathfrak{T}$ is called connected if for any distinct vertex $x$, $y$ on $\mathfrak{T}$, there exists a finite sequence of vertices $\{x_i\}_{i=0}^k$, $k\in\mathbb{Z}^+:=\{1,2,\ldots\}$ such that $x=x_0\sim x_1\sim\cdots\sim x_k=y$.
The geodesic path from the vertex $x$ to $y$ is denoted by $[x,y]$.
For notational simplicity, we identify $\mathfrak{T}$ with its set of vertices, $[x, y]$ with the set of all vertices lying on this geodesic path, and endow $\mathfrak{T}$ with the standard graph distance $d$, counting the number of edges in the shortest path connecting two vertices, and the set of its vertices with the counting measure $\mu$.
Under this construction, $(\mathfrak{T},d,\mu)$ becomes a metric measure space.

Let $\mathfrak{T}$ be a locally finite tree (i.e., $3\leq\nu(x) < \infty$ for all vertices $x$), noting that the degree function $\nu: \mathfrak{T} \to \mathbb{Z}^+$ may be unbounded globally. Two key subclasses are: the homogeneous tree $\mathfrak{T}_b$, defined by $\nu(x)=b+1$ for all $x$ (with $b\in\mathbb{Z}^+$, $b\geq2$); and trees of $(a,b)$-bounded geometry (class $\Upsilon_{a,b}$), satisfying $3\leq a+1\leq \nu(x)\leq b+1<+\infty$ for all $x$.

Within the metric measure space $(\mathfrak{T}, d, \mu)$ defined above, we introduce the fractional Hardy--Littlewood maximal operators. For a parameter $\gamma > 0$ and a function $f: \mathfrak{T} \to \mathbb{R}$, the centered and uncentered fractional Hardy--Littlewood maximal operators are defined respectively by
\[
\mathcal{M}^{\gamma}f(x) := \sup_{r \in \mathbb{N}} \frac{1}{\mu(B(x, r))^{\gamma}} \int_{B(x, r)} |f(y)| \, d\mu(y), \qquad x \in \mathfrak{T},
\]
and
\[
\widetilde{\mathcal{M}}^{\gamma}f(x) := \sup_{B \ni x} \frac{1}{\mu(B)^{\gamma}} \int_{B} |f(y)| \, d\mu(y), \qquad x \in \mathfrak{T}.
\]
Here,  $\mathbb{N}:=\{0,1,2,\ldots\}$. $B(x, r) := \{ y \in \mathfrak{T} : d(x, y) \le r \}$ is the closed ball of radius $r \ge 0$ centered at $x$, and the supremum in $\widetilde{\mathcal{M}}^{\gamma}$ is taken over all such closed balls $B$ containing $x$.
If $\gamma=1$, then $\mathcal{M}^{\gamma}$ and $\widetilde{\mathcal{M}}^{\gamma}$ are reduced to the classical centered and uncentered Hardy-Littlewood maximal operators $\mathcal{M}$ and $\widetilde{\mathcal{M}}$ on $\mathfrak{T}$.

In order to study the mapping properties of maximal operators on locally finite trees, we now introduce the definitions of $\ell^p$ space, $\ell^{p,\infty}$ space on the locally finite tree $\mathfrak{T}$.
\begin{dfn}{\rm($\ell^p(\mathfrak{T})$ space and $\ell^{p,\infty}(\mathfrak{T})$ space)} For $0<p\leq+\infty$, $\ell^p(\mathfrak{T})$ is the set of all functions $f:\mathfrak{T}\to\mathbb{R}$ on the locally finite tree $\mathfrak{T}$ such that $\Vert f\Vert_p<+\infty$, where
\begin{equation*}
\Vert f\Vert_p:=
\begin{cases}
\Big(\sum\limits_{j\in\mathfrak{T}}\vert f(j)\vert^p\Big)^{1/p} & {0<p<+\infty,}\\
\displaystyle\sup_{j\in\mathfrak{T}}\vert f(j)\vert & {p=+\infty,}
\end{cases}
\end{equation*}
and $\ell^{p,\infty}(\mathfrak{T})$ is the set of all functions $f:\mathfrak{T}\to\mathbb{R}$ on the locally finite tree $\mathfrak{T}$ such that $\Vert f\Vert_{p,\infty}<+\infty$, where
\begin{equation*}
\Vert f\Vert_{p,\infty}:=
\begin{cases}
\displaystyle\sup_{\lambda>0}\lambda\mu(\{x\in\mathfrak{T}:\vert f(x)\vert>\lambda\})^{{1}/{p}} & {0<p<+\infty,}\\
\displaystyle\sup_{j\in\mathfrak{T}}\vert f(j)\vert & {p=+\infty.}
\end{cases}
\end{equation*}
\end{dfn}

The study of the Hardy--Littlewood maximal operators on graphs was initiated by  Kor\'anyi and Picardello \cite{KP} in connection with the boundary behavior of eigenfunctions of the Laplace operator on trees.
Subsequently, Cowling, Meda, and Setti \cite[Theorem 3.1]{CMS} established that on the homogeneous tree $\mathfrak{T}_b$, the Hardy-–Littlewood maximal operator $\mathcal{M}$ is bounded on $\ell^p(\mathfrak{T}_b)$ for every $1<p\leq+\infty$ and is of weak type $(1,1)$.
Later, Veca \cite[Theorem 5.1]{V} proved that the fractional maximal operator $\mathcal{M}^{\frac12}$ on  $\mathfrak{T}_b$ is of restricted weak type $(2,2)$. 
For $0<\gamma\leq 1$, Levi and Santagati \cite[Theorem 3.4]{LS} established that  $\mathcal{M}^\gamma$ on $\mathfrak{T}_b$ is of strong type $(p,q)$ under either of these conditions: (i) $1\leq p\leq q\leq+\infty$, $q>{1}/{\gamma}$ and $p<{1}/{(1-\gamma)}$; (ii) $p={1}/{(1-\gamma)}$ and $q=+\infty$.
Using the volume growth property $\mu(B(x,r))\approx b^r$, one can compare
uncentered averages with centered fractional averages of order $1/2$.
This yields a pointwise domination of the uncentered maximal operator
$\widetilde{\mathcal M}$ by $\mathcal M^{1/2}$.
As a consequence, the $\ell^p\to \ell^q$ boundedness
of $\widetilde{\mathcal M}$ follows from the corresponding boundedness of
$\mathcal M^{1/2}$ for the same ranges of $p$ and $q$.

Beyond homogeneous trees, the study of maximal operators extends to diverse and non-trivial classes of graphs. Soria and Tradacete \cite[Theorem 4.9]{ST1} established a connection between the weak type $(1,1)$ estimates of $\mathcal{M}$ on the infinite graphs and some geometrical indices associated to the graph.   
In \cite[Theorem 1.5]{NT}, Naor and Tao proved that $\mathcal{M}$ is of weak type $(1,1)$ and bounded on $\ell^{p}(T^k)$ for every $1<p\leq+\infty$, where $T^k$ denotes the infinite rooted $k$-ary tree, $k\geq2$, that is, the infinite rooted tree such that each vertex has $k$ children. 
Subsequently, Soria and Tradacete \cite[Theorem 4.1]{ST2} showed that the proof structure of \cite[Theorem 1.5]{NT} can be adapted to nonhomogeneous trees and, more generally, to graphs with a globally bounded degree. 
Although the condition is highly versatile, verifying it in specific cases often requires considerable combinatorial or computational effort.
Note that $T^k\in\Upsilon_{k-1,k}$. 
For a general tree $\mathfrak{T}$ with bounded geometry in $\Upsilon_{a,b}$, where $2\leq a\leq b<+\infty$, Levi et al. \cite{LMSV} proved that if $a\leq b\leq a^2$, then $\mathcal{M}$ is bounded on $\ell^{p}(\mathfrak{T})$ for $\log_a b<p\leq+\infty$, and this lower bound of $p$ is sharp.
If $b>a^2$, then there exists a counterexample in $\Upsilon_{a,b}$ such that $\mathcal{M}$ is bounded on $\ell^{p}$ if and only if $p=+\infty$ (see \cite[Proposition 3.3]{LMSV}). 
Regarding $\widetilde{\mathcal{M}}$, it is proved in \cite[Theorem 4.1]{LMSV} that there are trees in $\Upsilon_{a,b}$ $(2\leq a<b<+\infty)$ for which $\widetilde{\mathcal{M}}$ is bounded on $\ell^p$ if and only if $p=+\infty$.
Moreover, there exist examples of trees of unbounded geometry \cite[Section 3]{LMSV} such that $\sup_{x\in\mathfrak{T}}\nu(x)=+\infty$ and $\mathcal{M}$ is bounded on $\ell^p$ if and only if $p=+\infty$.

The maximal operators $\mathcal{M}$ and $\widetilde{\mathcal{M}}$, defined by averaging over balls, are known to be bounded on $\ell^{p}$ ($1<p\leq+\infty$) for trees of bounded geometry. It is natural to ask if there are some other ``geometric" maximal operators on locally finite trees with potentially unbounded geometry that exhibit stable $\ell^p$ estimates for some $1\leq p\leq+\infty$. 
In \cite{MS}, Meda and Santagati defined  the centered and uncentered triangular maximal operators $\mathcal{N}$ and $\widetilde{\mathcal{N}}$, and stated that, if $\mathfrak{T}$ is a locally finite tree with $3\leq \nu(x)<+\infty$ for every $x\in\mathfrak{T}$, then $\mathcal{N}$ is bounded on $\ell^p(\mathfrak{T})$ for every $1\leq p\leq+\infty$ and $\widetilde{\mathcal{N}}$ is bounded on $\ell^p(\mathfrak{T})$ for every $1< p\leq+\infty$.
However, $\widetilde{\mathcal{N}}$ fails to be of even weak type $(1,1)$ on the homogeneous tree $\mathfrak{T}_b$, $b \geq 2$ see \cite[Theorems 3.1 \& 4.3]{MS}.

The aforementioned results on triangular maximal operators naturally lead to the question of their fractional analogues. Motivated by this, we study the mapping properties of fractional triangular maximal operators on locally finite trees, with a particular emphasis on establishing sharp boundedness results. We begin by setting up the relevant notations and definitions.

A geodesic ray $\omega$ in $\mathfrak{T}$ is a one-sided sequence of vertices $\{x_n:n\in\mathbb{N}\}$ of $\mathfrak{T}$ such that $x_{n+1}$ is the neighbor of $x_n$ for every $n$, and $d(x_i,x_j)=\vert i-j\vert$ for all $i,j\in\mathbb{N}$.
We write $x\in\omega$ if $x=x_n$ for some $n\in\mathbb{N}$ and we say that $x$ lies on $\omega$. For a fixed geodesic ray $\omega$ and a given vertex $y$, denote by $y\wedge\omega$ the vertex on $\omega$ closest to $y$. In other ways,
\[
y\wedge\omega:=\arg\min\limits_{x\in\omega}d(y,x).
\]
Denote by $\omega_y$ the geodesic ray $[y,x_k]\cup[x_k,x_{k+1},\ldots]$ if $y\wedge\omega=x_k$. For another vertex $x\in\mathfrak{T}$, denote by $x\succeq y$ if $x\in\omega_y$. If $x\succeq y$ and $x\neq y$, then it is denoted by $x\succ y$. Besides, we say that $x$ is a predecessor of $y$, and $y$ is a successor of $x$.

For a fixed geodesic ray $\omega=\{x_n,n\in\mathbb{N}\}$, the height function $h_\omega:\mathfrak{T}\to\mathbb{Z}$ is the discrete analogue of the Busemann function in Riemannian geometry, defined by 
\begin{align}\label{hf}
h_{\omega}(x):=\lim_{n\to+\infty}(n-d(x,x_n)),\; x\in\mathfrak{T}.
\end{align}
$\mathfrak{T}$ can be decomposed by the levels sets of the height function $h_\omega$:
\[\mathfrak{T}=\bigcup_{j\in\mathbb{Z}}\mathfrak{H}_{j}^{\omega},\]
in which $\mathfrak{H}_{j}^{\omega}$ is the $\omega$-horocycle at the height $j\in\mathbb{Z}$, defined by 
\[\mathfrak{H}_{j}^{\omega}:=\{x\in\mathfrak{T}:h_{\omega}(x)=j\}.\]

For $x\in\mathfrak{H}_{j}^{\omega}$, the direct successors of $x$ are vertices that are $\nu(x)-1$ neighbors of $x$ belonging to $\mathfrak{H}_{j-1}^{\omega}$. Oppositely, the direct predecessor of $x$ is the one neighbor of $x$ belonging to $\mathfrak{H}_{j+1}^{\omega}$. We denote by $p^{1}_{\omega}(x)$ and $s_{\omega}^{1}(x)$ the direct predecessor and the set of the direct successors of $x$ respectively. 
Set $s_{\omega}^{0}(x)=p_{\omega}^{0}(x):=\{x\}$, $s_\omega(x):=s_{\omega}^{1}(x)$, $p_{\omega}^{2}(x):=p_{\omega}(p_{\omega}(x))$. For positive integer $k\geq2$, we define the $k^{\rm th}$ predecessor and the $k^{\rm th}$ successor of $x$ iteratively:
\[
p_{\omega}^{k}(x):=p_{\omega}(p_{\omega}^{k-1}(x)),\quad s_{\omega}^{k}(x):=\bigcup_{y\in s_{\omega}^{k-1}(x)}s_{\omega}(y).
\]
For two distinct vertices $x$ and $y$ in $\mathfrak{T}$, their meet, denoted $x \wedge y$, is the unique vertex $v\in\mathfrak{T}$ satisfying: 
(i) $v$ is a common predecessor of both $x$ and $y$, that is, there exist integers $m, n \geq 1$ such that $x$ is an $m$-th successor of $v$ and $y$ is an $n$-th successor of $v$, and 
(ii) any vertex that is a common predecessor of $x$ and $y$ must be a predecessor of $v$, see \cite[Chapter 3]{D} for more details. The triangle $ T_{\omega}(x,r)$ with vertex $x$ and height $r\in\mathbb{N}$, is then defined by
\[
T_{\omega}(x,r):=\bigcup_{j=0}^{r}s_{\omega}^{j}(x),                                    
\]
where a triangle with vertex $x$ and height 0 is defined as $\{x\}$.

Let $\mathfrak{T}$ be a locally finite tree. We fix a reference vertex $o\in\mathfrak{T}$ and a geodesic ray $\omega$ in $\mathfrak{T}$, which induces a horocyclic decomposition of $\mathfrak{T}$. For notational simplicity, we denote the height function relative to $\omega$ by $h$, the horocycles by $\mathfrak{H}_j$ and write $p^{k}(x)$, $s^{k}(x)$, $T(x,r)$ as shorthand for $p^{k}_{\omega}(x)$, $s_{\omega}^{k}(x)$, $T_{\omega}(x,r)$ respectively.

For any triangle $T$, denote $v( T):=\arg\max_{x\in T}h(x)$ by the vertex of $T$, $h(T):=\max\limits_{x,y\in T}(h(x)-h(y))$ by the height of $T$, and $\beta(T):=s^{h( T)}(v( T))$ by the base of $T$, respectively.

We now introduce the definitions of fractional triangular maximal operators on locally finite trees $\mathfrak{T}$.
\begin{dfn} Let $\mathfrak{T}$ be a locally finite tree and let $\omega$ be a fixed geodesic ray in $\mathfrak{T}$. For any function $f: \mathfrak{T} \to \mathbb{R}$ and any $\gamma > 0$, the centered and the uncentered fractional triangular maximal operators are defined by
\[
\mathcal{N}_{\omega}^{\gamma}f(x)
:=\sup_{r\in\mathbb{N}}\frac{1}{\mu( T(x,r))^{\gamma}}\int_{ T(x,r)}\vert f(y)\vert d\mu(y), \quad x\in\mathfrak{T},
\]
and
\[
\widetilde{\mathcal{N}}_{\omega}^{\gamma}f(x):=\sup_{ T\ni x}\frac{1}{\mu( T)^{\gamma}}\int_{ T}\vert f(y)\vert d\mu(y), \quad x\in\mathfrak{T},
\]
respectively, where $T$ is any triangle containing $x$ in $\mathfrak{T}$.
For $\gamma=1$, $\mathcal{N}_{\omega}^{\gamma}$ and $\widetilde{\mathcal{N}}_{\omega}^{\gamma}$ are then reduced to the centered and uncentered triangular maximal operators in \cite{MS}.
\end{dfn}
Different from $\mathcal{M}^{\gamma}$ and $\widetilde{\mathcal{M}}^{\gamma}$ on locally finite trees, $\mathcal{N}_{\omega}^{\gamma}$ and $\widetilde{\mathcal{N}}_{\omega}^{\gamma}$ depend on the choice of the geodesic ray $\omega$.
In the subsequent analysis, we shall use the notation $\mathcal{N}^{\gamma}$ and $\widetilde{\mathcal{N}}^{\gamma}$ for the following reasons: all examples considered are either invariant under different choices of the geodesic ray $\omega$ (such as the homogeneous trees $\mathfrak{T}_b$ and the semi-homogeneous tree $\mathfrak{D}_{a,b}$, see Section 2), or are provided with an explicitly specified $\omega$ (for instance, the tree $\mathfrak{F}_{a,b}$ belonging to the class $\Upsilon_{a,b}$ as defined in Section 2).

For $\gamma=1$, it is established that $\mathcal{N}$ is of strong type $(p,q)$ if $1\leq p\leq q\leq +\infty$ \cite[Theorem 3.1]{MS}, while $\widetilde{\mathcal{N}}$ is of strong type $(p,q)$ for $1<p\leq q\leq +\infty$ by Marcinkiewicz interpolation and $\ell^p$ inclusions \cite[Theorem 4.3]{MS}. The case $\gamma>1$ is straightforward, yielding analogous boundedness results for $\mathcal{N}^\gamma$ and $\widetilde{\mathcal{N}}^\gamma$. 
Therefore, the nontrivial and interesting range is $0<\gamma\leq 1$, which will be the focus of our subsequent analysis.

Following the above line of investigation, the paper is structured in the following way. In Section \ref{s2}, we establish pointwise relations among the maximal operators $\mathcal{M}^\gamma$, $\widetilde{\mathcal{M}}^\gamma$, $\mathcal{N}^\gamma$, and $\widetilde{\mathcal{N}}^\gamma$ on locally finite trees. A central tool is a geometric construction of overlapping triangles with total measure comparable to a ball (Proposition \ref{lem21}), which we use to effectively estimate ball measures and derive the operator inequalities. Using this, we establish pointwise comparisons between these triangular maximal operators and the centered/uncentered fractional Hardy--Littlewood maximal operators on homogeneous trees (see Proposition \ref{p1}).
For general locally finite trees, we identify sufficient geometric conditions under which similar pointwise domination holds (see Proposition \ref{p2}). We also provide concrete counterexamples of locally finite trees to show that one maximal function may grow asymptotically much faster than another one. The results are summarized in Tables \ref{T1} and \ref{T2}.

In Section \ref{s3},  we investigate the boundedness of $\mathcal{N}^{\gamma}$ and $\widetilde{\mathcal{N}}^{\gamma}$ from $\ell^p$ to $\ell^q$ on locally finite trees. A key finding is that the uncentered operator $\widetilde{\mathcal{N}}^{\gamma}$ exhibits a distinct limitation: it fails to be of weak type $(1, 1/\gamma)$ on homogeneous trees (Theorem \ref{thm2}(iii)), in contrast to what might be expected by analogy with the centered case. 
Beyond this, we completely characterize the strong type $(p,q)$ boundedness for both operators (Theorems \ref{thm1}(i) and \ref{thm2}(i)), proving the obtained conditions on the exponents are optimal. The corresponding unboundedness results are given in Theorems \ref{thm1}(ii) and \ref{thm2}(ii).

As established in Section \ref{s2}, for $0<\gamma\leq1$, it remains an open question whether the inequality $\mathcal{N}^{\gamma}f(x)\lesssim\widetilde{\mathcal{M}}^{\gamma}f(x)$ holds for every $x$ in locally finite tree $\mathfrak{T}$ and all functions $f$.
Nevertheless, in Section \ref{s5},  we prove that on certain specific locally finite trees, such as $\mathfrak{F}_{a,b}$, $\widetilde{\mathcal{M}}^\gamma$ is unbounded for a strictly wider range of $(p,q)$ pairs than $\widetilde{\mathcal{N}}^\gamma$.

We use the notation $A\lesssim B$, whenever there exists constant $C>0$ (which is independent of the main parameters involved) such that $A\leq CB$.
Similarly for $A\gtrsim B$. We shall also use the notation $A\approx B$ for $A\lesssim B$ and $A\gtrsim B$. 
In this case, we say that $A$ and $B$ are equivalent or comparable. 

\section{Pointwise comparisons among  $\mathcal{M}^\gamma$, $\widetilde{\mathcal{M}}^\gamma$, $\mathcal{N}^\gamma$, and $\widetilde{\mathcal{N}}^\gamma$  on locally finite trees}\label{s2}

The main finding of this section is that the pointwise domination relations among $\mathcal{M}^\gamma$, $\widetilde{\mathcal{M}}^\gamma$, $\mathcal{N}^\gamma$, and $\widetilde{\mathcal{N}}^\gamma$ depend critically on the underlying tree geometry. We prove that on homogeneous trees $\mathfrak{T}_b$, the full chain $\mathcal{N}^\gamma\lesssim\mathcal{M}^\gamma \lesssim\widetilde{\mathcal{N}}^\gamma \lesssim\widetilde{\mathcal{M}}^\gamma$ holds. For general locally finite trees, we establish partial results under specific assumptions and provide counterexamples showing that the full chain may fail. The analysis relies on a geometric comparison between families of triangles and balls.

For any triangle $ T$ in any locally finite tree $\mathfrak{T}$ with $3 \leq a+1 \leq \nu < +\infty$, there is relation between the measure of the triangle and its base.

\begin{lem}\label{BASE2}
Let $\mathfrak{T}$ be a locally finite tree with $3 \leq a+1 \leq \nu < +\infty$, and let $T$ be a triangle in $\mathfrak{T}$ of height $h$. Then

(i) $a^k\mu(p^{k}(\beta( T)))\leq\mu(\beta( T))$ for every $k$ in $\{0,\ldots,h\}$;

(ii) $\mu( T)\leq\frac{a}{a-1}\mu(\beta( T))$.
\end{lem}

The proof of this lemma follows by adapting the method of \cite[Lemma 2.1]{MS}.

We begin by proving a proposition stating that every closed ball in a locally finite tree $\mathfrak{T}$ can be represented as a union of triangles or a union of disjoint triangles and trapezoids (see \cite[Definition 2]{ATV}), but not the admissible trapezoids (see \cite[Chapter 3]{HS}). Furthermore, the volume of the ball is comparable to the sum of the volumes of these triangles.

\begin{prop}\label{lem21}
Suppose that $\mathfrak{T}$ is a locally finite tree with $3\leq a+1\leq\nu<+\infty$.
Then for every $x\in\mathfrak{T}$ and $r\in\mathbb{N}$, we have that
\begin{align}\label{BUT}
B(x,r)=\bigcup\limits_{j=0}^{r}T(p^{j}(x),r-j)  
\end{align}
and 
\begin{align} \label{BUT2}
B(x,r)=p^r(x)\cup\bigcup_{j=0}^{r-1}\Big(\beta\left(T(p^{j}(x),r-j-1)\right)\cup
\beta\left(T(p^{j}(x),r-j)\right)\Big)
\end{align}
for $r\geq2$. Furthermore, for $r\in\mathbb{N}$,  it follows that
\begin{align}\label{MUBUT}
\frac{a-1}{a}\sum_{j=0}^{r}\mu(T(p^{j}(x),r-j))
\leq\mu(B(x,r))
\leq\sum_{j=0}^{r}\mu(T(p^{j}(x),r-j)).
\end{align}
\end{prop}

\begin{proof}
We firstly prove \eqref{BUT}.
Let $y\in T(p^{j_1}(x),r-j_{1})$ for some $0\leq j_{1}\leq r$ and $x\in\mathfrak{T}$.
By the definition of $T(p^{j_1}(x),r-j_{1})$, we have $d\big(p^{j_1}(x),y\big)\leq r-j_1$, and then
\[ 
d(x,y) \leq d\big(x,p^{j_1}(x)\big) + d\big(p^{j_1}(x),y\big)\le j_1 + (r-j_1) = r,
\]
implying $y\in B(x,r)$.

Conversely, for $y\in B(x,r)$, let $[y,x]$ be the unique geodesic path from $y$ to $x$.
Since $x\in\omega_x\bigcap[y,x]$, it follows that $\omega_x\bigcap[y,x]\neq\emptyset$.

We claim that there exists an integer $j_2\in[0,r]$ such that $y\wedge\omega_{x}=p^{j_2}(x)$.
In fact, for any vertex $z\in\omega_x$ one has $h(z)\ge h(x)$.  
Recall that for a geodesic ray $\omega=\{x_n\}_{n\in\mathbb{N}}$, the height
function is defined  in \eqref{hf}.
Since $z\in \omega_x$, the point \(y\) lies on the geodesic ray starting at
$x$ and pointing toward the limit of the sequence $\{x_n\}$. Therefore,
there exists an integer $N$ such that for all $n\geq N$, the three vertices
$x$, $z$, $x_n$ lie on the same geodesic and satisfy
\[
d(z,x_n)=d(x,x_n)-d(x,z), \qquad n\ge N.
\]
For such $n$ we compute
\[
n-d(z,x_n)
  = n - \bigl(d(x,x_n)-d(x,z)\bigr)
  = \bigl(n-d(x,x_n)\bigr) + d(x,z).
\]
Taking the limit as $n\to\infty$, and using the definition of $h$,
we obtain
\[
h(z)
= \lim_{n\to\infty}\bigl(n-d(z,x_n)\bigr)
 = \lim_{n\to\infty}\bigl(n-d(x,x_n)\bigr) + d(x,z)
 = h(x)+d(x,z).
\]
Since \(d(x,z)\geq 0\), this shows that
\[
h(z)\geq h(x),
\]
with equality if and only if $x=z$. 
Let $ z = y \wedge \omega_x $. Then $ z = p^{j_2}(x) $ for some nonnegative integer $ j_2 $. If, to the contrary, $ j_2 > r $, then
\[
d(y,x)=d(y,z)+d(z,x) > r,
\]
contradicting $ y \in B(x,r) $. Hence $ j_2 \leq r $, which proves the claim.

If $j_2=0$, then $x\succeq y$ and $y\in T(x,r)$. If $1\leq j_2\leq r$, then 
\[ 
d\big(p^{j_2}(x),y\big)=d(x,y)-j_2\leq r-j_2,
\]
hence $y\in T\big(p^{j_2}(x),r-j_2\big)$. 

For \eqref{BUT2}, if $r=2$, then we have that 
\begin{align*}
B(x,2)
&=\bigcup_{j=0}^{2}\Big(T(p^{j}(x),r-j)\Big)\\
&=p^2(x)\cup\Big(\beta(T(x,1))\cup\beta(T(x,2))\Big)
\cup\Big(\beta(T(p(x),0))\cup\beta(T(p(x),1))\Big)\\
&=p^2(x)\cup\bigcup_{j=0}^{1}\Big(\beta\left(T(p^{j}(x),2-j-1)\right)\cup
\beta\left(T(p^{j}(x),2-j)\right)\Big).
\end{align*}
If $r\geq3$, then it follow that  
\begin{align*}
B(x,r)
&=T(p^r(x),0)
\cup T(x,r)
\cup\Big(\bigcup_{j=1}^{r-1}
T(p^{j}(x),r-j)\setminus\Big(\bigcup_{i=0}^{j-1}T(p^{i}(x),r-i)\Big)\Big)\\
&=p^{r}(x)\cup T(x,r)
\cup\Big(
\bigcup_{j=1}^{r-1}\Big(T(p^{j}(x),r-j)\setminus T(p^{j-1}(x),r-j-1)\Big)\Big)\\
&=p^{r}(x)\cup T(p^{r-1}(x),1)
\cup\Big(\bigcup_{j=0}^{r-2}\Big(T(p^{j}(x),r-j)\setminus T(p^{j}(x),r-j-2)\Big)\Big)\\
&=p^{r}(x)
\cup\Big(\bigcup_{j=0}^{r-1}\Big(\beta\left(T(p^{j}(x),r-j-1)\right)
\cup\beta\left(T(p^{j}(x),r-j)\right)\Big)\Big).
\end{align*}

Now for \eqref{MUBUT}, for every $r\in\mathbb{N}$, we have that 
\begin{align*}
\mu(B(x,r))
=\mu\Big(\bigcup\limits_{j=0}^{r}T(p^{j}(x),r-j)\Big)
\leq\sum_{j=0}^{r}\mu\big(T(p^{j}(x),r-j)\big).
\end{align*}

We now verify the formula for $\mu(B(x, r))$ by distinguishing cases.

Case $r=0$: Trivially, $\mu(B(x,0)) = \mu(T(x,0)) = 1$.

Case $r=1$: By definition, $\mu(B(x,1)) = \mu(T(x,1)) + \mu(T(p(x),0))$.

Case $r=2$: We find
\[
\mu(B(x,2))
=\sum_{j=0}^{2}\Big(\mu(T(p^{j}(x),r-j))\Big)-1
\geq\frac{1}{2}\sum_{j=0}^{2}\mu(T(p^{j}(x),r-j)).
\]

Case $r \geq 3$: In this general case, we have
\begin{align*}
\mu(B(x,r))
&=1
+\sum_{j=0}^{r-1}\Big(\mu\left(\beta\left(T(p^{j}(x),r-j-1)\right)\right)
+\mu\left(\beta\left(T(p^{j}(x),r-j)\right)\right)\Big)\\
&\geq1
+\sum_{j=0}^{r-1}\mu\left(\beta\left(T(p^{j}(x),r-j)\right)\right).
\end{align*}
From Lemma \ref{BASE2} (ii), we have that
\[
\mu\left(\beta\left(T(p^{j}(x),r-j)\right)\right)
\geq
\frac{a-1}{a}\mu\left(T(p^{j}(x),r-j)\right).
\]
Therefore, we have that
\begin{align*}
\mu(B(x,r))
\geq
\frac{a-1}{a}\Big(1+\sum_{j=0}^{r-1}\mu\left(T(p^{j}(x),r-j)\right)\Big)
=
\frac{a-1}{a}\sum_{j=0}^{r}\mu(T(p^{j}(x),r-j)).
\end{align*}
This completes the proof of \eqref{MUBUT}.
\end{proof}

\begin{rem}
\cite[Lemma 2.1]{LMSV} provides a volume formula for spheres in trees of bounded geometry, which can be applied to estimate the volume of balls. On locally finite trees, however, the volume of spheres and balls are not necessarily comparable. For this reason, \eqref{MUBUT} offers an alternative method to evaluate the volume of balls on locally finite trees. 
\end{rem}

Now we give the pointwise comparisons among operators $\mathcal{M}^\gamma$, $\widetilde{\mathcal{M}}^\gamma$,  $\mathcal{N}^\gamma$ and $\widetilde{\mathcal{N}}^\gamma$ on homogeneous trees.

\begin{prop}\label{p1}
	Let $\gamma > 0$ and let $\mathfrak{T}_b$ be a homogeneous tree with $b \geq 2$. Then for every function $f$ on $\mathfrak{T}_b$ and every $x \in \mathfrak{T}_b$,  we have
	\begin{align}\label{C1}
		\mathcal{N}^{\gamma}f(x) \lesssim \mathcal{M}^{\gamma}f(x) \lesssim \widetilde{\mathcal{N}}^{\gamma}f(x) \lesssim \widetilde{\mathcal{M}}^{\gamma}f(x).
	\end{align}
\end{prop}

\begin{proof}
For $r\geq1$, $x\in\mathfrak{T}_b$, we have
\begin{align*}
\mu(B(x,r))=\frac{b^{r+1}+b^r-2}{b-1},\qquad\mu( T(x,r))=\frac{b^{r+1}-1}{b-1},
\end{align*}
and $\mu(B(x,0))=\mu( T(x,0))=1$. Then we can obtain that 
\[
\sup_{r\geq0}\frac{\mu(B(x,r))}{\mu( T(x,r))}
=\sup_{r\geq1}\frac{b+1-\frac{2}{b^r}}{b-\frac{1}{b^r}}=1+\frac{1}{b}.
\]
Moreover, from the definition of $\mathfrak{T}_b$, for any ball $B_r$ and triangle $T_r$ containing $x$ with radius and height $r$, $\mu(B_r)=\mu(B(x,r))\approx\mu(T(x,r))=\mu(T_r)\approx b^r$, $x\in\mathfrak{T}_b$ and $r\in\mathbb{N}$.

Then for any $r\geq0$ and $f$ on $\mathfrak{T}_b$, we have that 
\begin{align*}
\frac{1}{\mu( T(x,r))^{\gamma}}\int_{ T(x,r)}\vert f(y)\vert d\mu(y)
\lesssim\frac{1}{\mu(B(x,r))^\gamma}\int_{B(x,r)}\vert f(y)\vert d\mu(y). 
\end{align*}
By taking the supremum over all $r\geq0$, one has that $\mathcal{N}^{\gamma}f(x)
\lesssim\mathcal{M}^{\gamma}f(x)$ and similarly we obtain  $\widetilde{\mathcal{N}}^{\gamma}f(x)
\lesssim\widetilde{\mathcal{M}}^{\gamma}f(x)$ holds for every $x\in\mathfrak{T}_b$ and all $f$ on $\mathfrak{T}_b$.

Denote by $\lfloor x \rfloor = \max \{ n \in \mathbb{Z}: n \le x \}$. Then it follows from \eqref{BUT} that
\begin{align*}
\frac{1}{\mu(B(x,r))^{\gamma}}&\int_{B(x,r)}\vert f(y)\vert d\mu(y)
=\frac{1}{\mu(B(x,r))^{\gamma}}
\int_{\bigcup\limits_{j=0}^{r}T(p^{j}(x),r-j)}\vert f(y)\vert d\mu(y)\\
&\hspace{-1cm}
\leq\sum_{j=0}^{\lfloor r/2\rfloor}
\frac{\mu(T(p^{j}(x),r-j))^{\gamma}}{\mu(B(x,r))^{\gamma}}
\frac{1}{\mu(T(p^{j}(x),r-j))^{\gamma}}\int_{T(p^{j}(x),r-j)}\vert f(y)\vert d\mu(y)\\
&\hspace{-0.7cm}
+\sum_{j=\lfloor r/2\rfloor+1}^{r}
\frac{\mu(T(p^{j}(x),j))^{\gamma}}{\mu(B(x,r))^{\gamma}}
\frac{1}{\mu(T(p^{j}(x),j))^{\gamma}}\int_{T(p^{j}(x),j)}\vert f(y)\vert d\mu(y)\\
&\hspace{-1cm}
\leq\Big(\sum_{j=0}^{\lfloor r/2\rfloor}
\frac{\mu(T(p^{j}(x),r-j))^{\gamma}}{\mu(B(x,r))^{\gamma}}
+\sum_{j=\lfloor r/2\rfloor+1}^{r}
\frac{\mu(T(p^{j}(x),j))^{\gamma}}{\mu(B(x,r))^{\gamma}}\Big)
\widetilde{\mathcal{N}}^{\gamma}f(x)\\
&\hspace{-1cm}
\approx\Big(\sum_{j=0}^{\lfloor r/2\rfloor}b^{-j\gamma}
+\sum_{j=\lfloor r/2\rfloor+1}^{r}b^{(j-r)\gamma}\Big)
\widetilde{\mathcal{N}}^{\gamma}f(x)\\
&\hspace{-1cm}
\leq2\Big(\sum_{j=0}^{+\infty}b^{-j\gamma}\Big)
\widetilde{\mathcal{N}}^{\gamma}f(x)
=\frac{2b^\gamma}{b^\gamma-1}\widetilde{\mathcal{N}}^{\gamma}f(x).
\end{align*}
Therefore, \eqref{C1} holds. This completes the proof.
\end{proof}
\vspace{-0.5em}
Based on the conclusions of Proposition \ref{p1}, a natural question arises: whether certain pairs among the four operators in \eqref{C1} are equivalent on homogeneous trees. The answer is negative, as demonstrated by the following result.

\begin{prop}\label{phc}
Suppose that $\mathfrak{T}_b$ $(b\geq2)$ is a homogeneous tree. Let $\gamma>0$, $n\in\mathbb{Z}^+$, $x_0\in\mathfrak{T}_b$ be the reference vertex.

(i) For every vertex $x^\ast\in\mathfrak{T}_b$, define the  function
\begin{align*}
f_{x^\ast}(x):=\sum_{k=1}^{+\infty}\chi_{p^{k}(x^\ast)}(x),\; x\in\mathfrak{T}_b.
\end{align*}
Then the following relationships fail:
\begin{align}
\mathcal{M}^{\gamma}f_{x^\ast}(x)&\lesssim \mathcal{N}^{\gamma}f_{x^\ast}(x), \label{phc11}\\
\widetilde{\mathcal{N}}^{\gamma}f_{x^\ast}(x)
&\lesssim \mathcal{N}^{\gamma}f_{x^\ast}(x), \label{phc13}\\
\widetilde{\mathcal{M}}^{\gamma}f_{x^\ast}(x)
&\lesssim \mathcal{N}^{\gamma}f_{x^\ast}(x).\label{phc12}
\end{align}

(ii) Let $x_n$ be a vertex belonging to $s^{n}(x_0)$. The Kronecker delta function $\delta_{x_n} : \mathfrak{T}_b \to \{0,1\}$ is defined as
\begin{equation*}
\delta_{x_n}(x)=
\begin{cases}
1 & {x=x_n,}\\
0 & {x\neq x_n.}
\end{cases}
\end{equation*}
Then
\begin{align}\label{p22}
\lim_{n\to+\infty}
\frac{\widetilde{\mathcal{M}}^{\gamma}\delta_{x_n}(x_0)}
{\mathcal{M}^{\gamma}\delta_{x_n}(x_0)}
=\frac{\widetilde{\mathcal{M}}^{\gamma}\delta_{x_n}(x_0)}
{\widetilde{\mathcal{N}}^{\gamma}\delta_{x_n}(x_0)}
=+\infty.
\end{align}

(iii) Let $x_n^\ast\in s^n(p^n(x_0))$ such that $x_n^\ast\wedge x_0=p^{n}(x_0)$ and $d(x_n^\ast,x_0)=2n$. Then
\begin{align}\label{p23}
\lim_{n\to+\infty}
\frac{\widetilde{\mathcal{N}}^{\gamma}\delta_{x_n^\ast}(x_0)}
{\mathcal{M}^{\gamma}\delta_{x_n^\ast}(x_0)}
=+\infty.
\end{align}

\end{prop}

\begin{proof}
(i) For every $r\geq0$, $x\in T(x^\ast,r)$, we have 
\[\int_{ T(x^\ast,r)} f_{x^\ast}(x) d\mu(x)=0,\]
which implies that $\mathcal{N}^{\gamma}f_{x^\ast}(x)\equiv0$. Nevertheless, since
\begin{align*}
\mathcal{M}^{\gamma}f_{x^{\ast}}(x)
&=\sup_{r\geq0}\frac{1}{\mu(B(x,r))^{\gamma}}\int_{B(x,r)}\vert f_{x^\ast}(y)\vert d\mu(y)
\geq\sup_{k\in\mathbb{Z}^+}\frac{k}{\mu(B(x,d(x,p^{k}(x^\ast))))^{\gamma}}
\\
&\approx\sup_{k\in\mathbb{Z}^+}\frac{k}{b^{d(x,p^{k}(x^\ast))\gamma}}
=\sup_{k\in\mathbb{Z}^+}\frac{k}{b^{(d(x,x^\ast)+d(x^\ast,p^{k}(x^\ast)))\gamma}}\\
&=\frac{1}{b^{d(x,x^\ast)\gamma}}\sup_{k\in\mathbb{Z}^+}\frac{k}{b^{k\gamma}}
=\frac{1}{b^{(d(x,x^\ast)+1)\gamma}},
\end{align*}
there does not exist a constant $C>0$ such that \eqref{phc11} holds.
Since $\widetilde{\mathcal{M}}^{\gamma}f_{x^\ast}(x)
\gtrsim\widetilde{\mathcal{N}}^{\gamma}f_{x^\ast}(x)
\gtrsim\mathcal{M}^{\gamma}f_{x^\ast}(x)$, it follows that neither \eqref{phc12} nor \eqref{phc13} holds.

(ii) Without loss of generality, we can assume that $n\in\mathbb{Z}^+$ is even. Then we have that 
\begin{align*}
\frac{\widetilde{\mathcal{M}}^{\gamma}\delta_{x_n}(x_0)}
{\widetilde{\mathcal{N}}^{\gamma}\delta_{x_n}(x_0)}
&=\frac{\displaystyle\frac{1}{\mu(B(p^{n/2}(x_n),n/2))^\gamma}
\int_{B(p^{n/2}(x_n),n/2)}\delta_{x_n}(y)d\mu(y)}
{\displaystyle\frac{1}{\mu( (T(x_0,n))^\gamma}\int_{ T(x_0,n)}\delta_{x_n}(y)d\mu(y)}\\
&=\Big(\frac{\mu( T(x_0,n))}{\mu(B(p^{n/2}(x_n),n/2))}\Big)^\gamma\\
&\approx b^{\frac{n\gamma}{2}}\to+\infty
\end{align*}
as $n\to+\infty$. Since $\mathcal{M}^{\gamma}\delta_{x_n}(x)
\lesssim\widetilde{\mathcal{N}}^{\gamma}\delta_{x_n}(x)$ for any $x\in\mathfrak{T}_b$, it follows that 
\[\frac{\widetilde{\mathcal{M}}^{\gamma}\delta_{x_n}(x_0)}
{\mathcal{M}^{\gamma}\delta_{x_n}(x_0)}
\gtrsim\frac{\widetilde{\mathcal{M}}^{\gamma}\delta_{x_n}(x_0)}
{\widetilde{\mathcal{N}}^{\gamma}\delta_{x_n}(x_0)},\]
which implies that \eqref{p22} holds.

(iii) Indeed, it follows from the following calculation:
\begin{align*}
\frac{\widetilde{\mathcal{N}}^{\gamma}\delta_{x_n^\ast}(x_0)}
{\mathcal{M}^{\gamma}\delta_{x_n^\ast}(x_0)}
&=\frac{\displaystyle\frac{1}{\mu( T(p^n(x_0),n))^\gamma}
\int_{ T(p^n(x_0),n)}\delta_{x_n^\ast}(y)d\mu(y)}
{\displaystyle\frac{1}{\mu(B(x_0,2n))^\gamma}\int_{B(x_0,2n)}\delta_{x_n^\ast}(y)d\mu(y)}\\
&=\left(\frac{\mu(B(x_0,2n))}{\mu( T(p^n(x_0),n))}\right)^\gamma\\
&\approx b^{n\gamma}\to+\infty
\end{align*}
as $n\to+\infty$, which implies that \eqref{p23} holds.
\end{proof}
\vspace{-0.5em}
Apart from homogeneous trees, there exist further examples of trees that satisfy \eqref{C1} as well. By checking the proof of Proposition \ref{p1}, we have the following result. 
\begin{cor}
Let $\mathfrak{T}$ be a locally finite tree with $3\leq a+1\leq\nu<+\infty$, $\gamma>0$. Assume that there exists a positive constant $\tau>1$ such that
\[
\mu(T(x,r))\approx\tau^r,\quad \forall x\in\mathfrak{T},\; r\in\mathbb{N}.
\]
Then \eqref{C1} holds for every $x\in\mathfrak{T}$.
\end{cor}

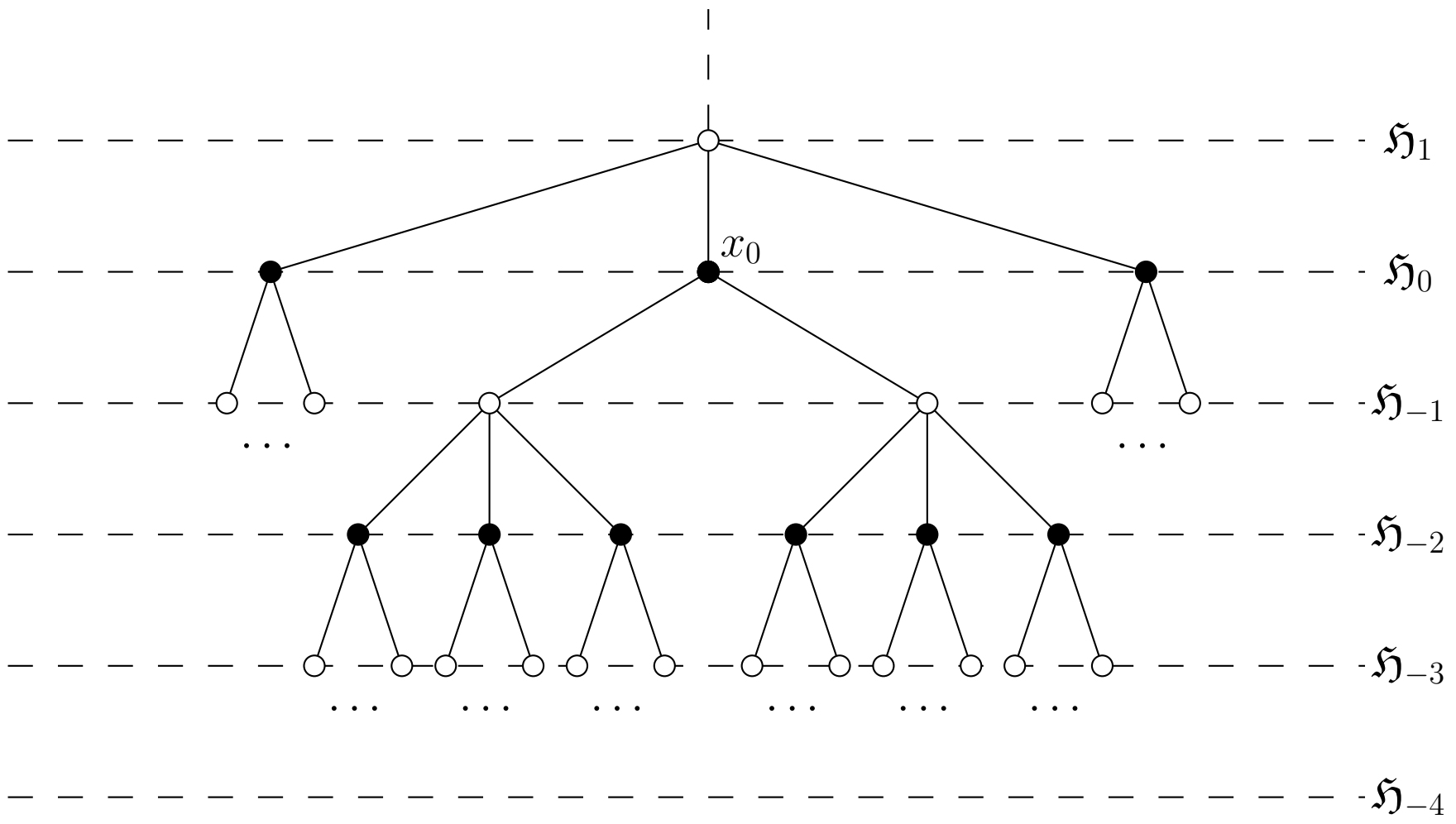
\begin{figure}[!htbp]
  \centering
\tikzstyle{every node}=[scale=0.5]
\begin{tikzpicture}[scale=0.35]
	\begin{pgfonlayer}{nodelayer}
		\node [style=Black] (8) at (0, 3) {};
		\node [style=Black] (9) at (0, 3) {};
		\node [style=Black] (23) at (0, 3) {};
		\node [style=new style 0] (24) at (5, 0) {};
		\node [style=Black] (26) at (2, -3) {};
		\node [style=new style 0] (28) at (3, -6) {};
		\node [style=new style 0] (29) at (1, -6) {};
		\node [style=new style 0] (33) at (7, -6) {};
		\node [style=Black] (35) at (8, -3) {};
		\node [style=new style 0] (36) at (9, -6) {};
		\node [style=new style 0] (37) at (0, 6) {};
		\node [style=Black] (38) at (10, 3) {};
		\node [style=new style 0] (39) at (-1, -6) {};
		\node [style=new style 0] (40) at (-3, -6) {};
		\node [style=Black] (42) at (-2, -3) {};
		\node [style=new style 0] (47) at (-7, -6) {};
		\node [style=new style 0] (49) at (-9, -6) {};
		\node [style=Black] (50) at (-8, -3) {};
		\node [style=new style 0] (51) at (-5, 0) {};
		\node [style=Black] (52) at (-10, 3) {};
		\node [style=none] (53) at (0, 9) {};
		\node [style=none] (62) at (-16, 6) {};
		\node [style=none] (63) at (15, 6) {};
		\node [style=none] (64) at (-16, 3) {};
		\node [style=none] (65) at (15, 3) {};
		\node [style=none] (66) at (-16, 0) {};
		\node [style=none] (67) at (15, 0) {};
		\node [style=none] (68) at (-16, -3) {};
		\node [style=none] (69) at (15, -3) {};
		\node [style=new style 0] (70) at (-11, 0) {};
		\node [style=new style 0] (72) at (-9, 0) {};
		\node [style=new style 0] (73) at (9, 0) {};
		\node [style=new style 0] (75) at (11, 0) {};
		\node [style=none] (76) at (-16, -6) {};
		\node [style=none] (77) at (15, -6) {};
		\node [style=none] (78) at (0.75, 3.5) {\Huge$x_0$};
		\node [style=none] (79) at (16, 3) {};
		\node [style=none] (80) at (16, 3) {\Huge$\mathfrak{H}_0$};
		\node [style=none] (81) at (16, 6) {\Huge$\mathfrak{H}_1$};
		\node [style=none] (83) at (16, 0) {\Huge$\mathfrak{H}_{-1}$};
		\node [style=none] (84) at (16, -3) {\Huge$\mathfrak{H}_{-2}$};
		\node [style=none] (85) at (16, -6) {\Huge$\mathfrak{H}_{-3}$};
		\node [style=none] (87) at (-10, -1) {\Huge$\cdots$};
		\node [style=none] (88) at (10, -1) {\Huge$\cdots$};
		\node [style=none] (89) at (-8, -7) {\Huge$\cdots$};
		\node [style=none] (91) at (-2, -7) {\Huge$\cdots$};
		\node [style=none] (94) at (2, -7) {\Huge$\cdots$};
		\node [style=none] (95) at (5, -7) {\Huge$\cdots$};
		\node [style=none] (96) at (8, -7) {\Huge$\cdots$};
		\node [style=Black] (101) at (5, -3) {};
		\node [style=Black] (102) at (-5, -3) {};
		\node [style=new style 0] (103) at (-6, -6) {};
		\node [style=new style 0] (104) at (-4, -6) {};
		\node [style=new style 0] (107) at (4, -6) {};
		\node [style=new style 0] (108) at (6, -6) {};
		\node [style=none] (109) at (-16, -9) {};
		\node [style=none] (110) at (15, -9) {};
		\node [style=none] (111) at (16, -9) {\Huge$\mathfrak{H}_{-4}$};
		\node [style=none] (112) at (-5, -7) {\Huge$\cdots$};
	\end{pgfonlayer}
	\begin{pgfonlayer}{edgelayer}
		\draw (51) to (50);
		\draw (50) to (49);
		\draw (50) to (47);
		\draw (51) to (42);
		\draw (42) to (40);
		\draw (42) to (39);
		\draw (51) to (23);
		\draw (37) to (23);
		\draw (37) to (38);
		\draw (37) to (52);
		\draw [style=new edge style 0] (62.center) to (63.center);
		\draw [style=new edge style 0] (68.center) to (69.center);
		\draw [style=new edge style 0] (66.center) to (67.center);
		\draw [style=new edge style 0] (64.center) to (65.center);
		\draw [style=new edge style 0] (76.center) to (77.center);
		\draw [style=new edge style 0] (37) to (53.center);
		\draw (29) to (26);
		\draw (28) to (26);
		\draw (33) to (35);
		\draw (36) to (35);
		\draw (35) to (24);
		\draw (26) to (24);
		\draw (24) to (23);
		\draw (70) to (52);
		\draw (72) to (52);
		\draw (75) to (38);
		\draw (73) to (38);
		\draw (102) to (51);
		\draw (101) to (24);
		\draw (102) to (103);
		\draw (102) to (104);
		\draw (101) to (107);
		\draw (101) to (108);
		\draw [style=new edge style 0] (109.center) to (110.center);
	\end{pgfonlayer}
\end{tikzpicture}
  \caption{An example of $\mathfrak{D}_{a,b}$ when $a=2$, $b=3$}\label{Dab}
\end{figure}

A representative example is the semi-homogeneous tree $\mathfrak{D}_{a,b}$ (see \cite[p. 3966]{LMSV}) characterized by two degrees $a+1$ and $b+1$. Each vertex has either $a+1$ or $b+1$ neighbors and no two adjacent vertices share the same degree (see Figure \ref{Dab}). It is straightforward to check that 
\[
\mu(T(x,r))\approx(\sqrt{ab})^r,\quad \forall x\in\mathfrak{D}_{a,b},\; r\in\mathbb{N}.
\]

We establish a more general sufficient condition such that \eqref{C1} holds. 

\begin{prop}\label{p2}
Let $\gamma>0$. Suppose $\mathfrak{T}$ is a locally finite tree with $3\leq a+1\leq\nu<+\infty$. If there exists a positive constant $C_0\geq1$ such that 
\begin{align}\label{GMB}
\sup_{\substack{x \in \mathfrak{T} \\ r \geq 1, j \geq 1}}\frac{\mu(T(p^{j}(x),r))}{\mu(T(x,r))}\leq C_0<+\infty,
\end{align}
then \eqref{C1} holds for every $x\in\mathfrak{T}$.
\end{prop}

\begin{proof}
First, we claim that if \eqref{GMB} holds, then $\mu(B(x,r))\approx\mu(T(x,r))$ for every $x\in\mathfrak{T}$ and $r\in\mathbb{N}$. In fact, Lemma \ref{BASE2} and Proposition \ref{lem21} imply that 
\begin{align*}
\frac{\mu(B(x,r))}{\mu(T(x,r))}
&\approx\sum_{j=0}^{r}\frac{\mu(T(p^j(x),r-j))}{\mu(T(x,r))}\\
&=\sum_{j=0}^{r}\frac{\mu(T(p^j(x),r-j))}{\mu(T(p^j(x),r))}
\frac{\mu(T(p^j(x),r))}{\mu(T(x,r))}\\
&\leq C_0\sum_{j=0}^{r}\Big(\frac{1}{a}\Big)^j<\frac{C_0a}{a-1}.
\end{align*}
Since $T(x,r)\subset B(x,r)$ for every $x\in\mathfrak{T}$ and $r\in\mathbb{N}$, it follows that 
\[
\mu(T(x,r))
\leq\mu(B(x,r))
\leq\frac{C_0a}{a-1}\mu(T(x,r))
\]
for every $x\in\mathfrak{T}$ and $r\in\mathbb{N}$. Hence, we have that 
\begin{align*}
\frac{1}{\mu( T(x,r))^{\gamma}}\int_{ T(x,r)}\vert f(y)\vert d\mu(y)
&\lesssim\frac{1}{\mu(B(x,r))^\gamma}\int_{B(x,r)}\vert f(y)\vert d\mu(y)
\leq\mathcal{M}^{\gamma}f(y),
\end{align*}
and we obtain that $\mathcal{N}^{\gamma}f(x)
\lesssim\mathcal{M}^{\gamma}f(x)$. Similarly we can prove  $\widetilde{\mathcal{N}}^{\gamma}f(x)
\lesssim\widetilde{\mathcal{M}}^{\gamma}f(x)$ holds for any $x\in\mathfrak{T}$ and all $f$ on $\mathfrak{T}$. Furthermore, since 
\[
\bigcup_{j=\lfloor r/2\rfloor+1}^{r}T(p^j(x),r-j)\subset T(p^r(x),r),
\]
 the following estimates holds for any $x\in\mathfrak{T}$ and $r\geq0$:
\begin{align*}
\frac{1}{\mu(B(x,r))^{\gamma}}&\int_{B(x,r)}\vert f(y)\vert d\mu(y)
=\frac{1}{\mu(B(x,r))^{\gamma}}
\int_{\bigcup\limits_{j=0}^{r}T(p^{j}(x),r-j)}\vert f(y)\vert d\mu(y)\\
&\hspace{-1cm}
\leq\sum_{j=0}^{\lfloor r/2\rfloor}
\frac{\mu(T(p^{j}(x),r-j))^{\gamma}}{\mu(B(x,r))^{\gamma}}
\frac{1}{\mu(T(p^{j}(x),r-j))^{\gamma}}\int_{T(p^{j}(x),r-j)}\vert f(y)\vert d\mu(y)\\
&\hspace{-0.7cm}
+\frac{\mu(T(p^{r}(x),r))^{\gamma}}{\mu(B(x,r))^{\gamma}}
\frac{1}{\mu(T(p^{r}(x),r))^{\gamma}}\int_{T(p^{r}(x),r)}\vert f(y)\vert d\mu(y)\\
&\hspace{-1cm}
\approx\left(\sum_{j=0}^{\lfloor r/2\rfloor}\Big(
\frac{\mu(T(p^{j}(x),r-j))^{\gamma}}{\mu(T(x,r))^{\gamma}}\Big)
+\frac{\mu(T(p^{r}(x),r))^{\gamma}}{\mu(T(x,r))^{\gamma}}\right)
\widetilde{\mathcal{N}}^{\gamma}f(x)\\
&\hspace{-1cm}
\leq\left(C_0^\gamma+\sum_{j=0}^{\lfloor r/2\rfloor}
\frac{\mu(T(p^{j}(x),r-j))^{\gamma}}{\mu(T(p^{j}(x),r))^{\gamma}}
\frac{\mu(T(p^{j}(x),r))^{\gamma}}{\mu(T(x,r))^{\gamma}}\right)
\widetilde{\mathcal{N}}^{\gamma}f(x)\\
&\hspace{-1cm}
\leq\Big(C_0^\gamma+\sum_{j=0}^{\lfloor r/2\rfloor}\frac{C_0^\gamma}{a^{\gamma j}}\Big)
\widetilde{\mathcal{N}}^{\gamma}f(x)\\
&\hspace{-1cm}
\leq C_0^\gamma\Big(1+\frac{1}{1-a^{-\gamma}}\Big)\widetilde{\mathcal{N}}^{\gamma}f(x).
\end{align*}
We conclude the proof.
\end{proof}
\vspace{-0.5em}

From the proof of Proposition \ref{p2}, \eqref{GMB} implies that
\begin{align}\label{BTC}
\mu(B(x,r))\approx\mu(T(x,r))
\end{align}
for every $x\in\mathfrak{T}$ and $r\in\mathbb{N}$. It is still unknown whether \eqref{BTC} implies \eqref{C1}.
We state this open problem as the following conjecture.
\begin{con}
Suppose that $\gamma>0$, $\mathfrak{T}$ is a locally finite tree with $3\leq a+1\leq\nu<+\infty$. Assume that
\begin{align*}
\mu(B(x,r))\approx\mu( T(x,r))
\end{align*}
holds for every $x\in\mathfrak{T}$ and every $r\geq0$.
Then we have that
\begin{align*}
\mathcal{M}^{\gamma}f(x)\lesssim\widetilde{\mathcal{N}}^{\gamma}f(x)
\end{align*}
for every $x\in\mathfrak{T}$.
\end{con}

\begin{figure}[!htbp]
  \centering
\tikzstyle{every node}=[scale=0.5]
\begin{tikzpicture}[scale=0.23]
	\begin{pgfonlayer}{nodelayer}
		\node [style=none] (80) at (26, 1) {\Huge$\mathfrak{H}_0$};
		\node [style=none] (81) at (26, 5) {\Huge$\mathfrak{H}_1$};
		\node [style=none] (82) at (26, 9) {\Huge$\mathfrak{H}_2$};
		\node [style=none] (83) at (26, -3) {\Huge$\mathfrak{H}_{-1}$};
		\node [style=none] (84) at (26, -7) {\Huge$\mathfrak{H}_{-2}$};
		\node [style=none] (85) at (26, -11) {\Huge$\mathfrak{H}_{-3}$};
		\node [style=none] (87) at (22, 0) {\Huge$\cdots$};
		\node [style=none] (88) at (14, -8) {\Huge$\cdots$};
		\node [style=none] (89) at (-22, 0) {\Huge$\cdots$};
		\node [style=none] (90) at (-22, -8) {\Huge$\cdots$};
		\node [style=none] (91) at (-18, -8) {\Huge$\cdots$};
		\node [style=none] (92) at (-10, -8) {\Huge$\cdots$};
		\node [style=none] (93) at (-2, -8) {\Huge$\cdots$};
		\node [style=none] (94) at (2, -8) {\Huge$\cdots$};
		\node [style=none] (95) at (6, -8) {\Huge$\cdots$};
		\node [style=none] (96) at (10, -8) {\Huge$\cdots$};
		\node [style=none] (104) at (-6, -8) {\Huge$\cdots$};
		\node [style=none] (106) at (-14, -8) {\Huge$\cdots$};
		\node [style=new style 0] (112) at (-10, -3) {};
		\node [style=new style 0] (113) at (-11.5, -7) {};
		\node [style=new style 0] (114) at (-9.5, -7) {};
		\node [style=new style 0] (115) at (-13.5, -7) {};
		\node [style=new style 0] (116) at (-15.5, -7) {};
		\node [style=new style 0] (117) at (-14, -3) {};
		\node [style=new style 0] (118) at (-17.5, -7) {};
		\node [style=new style 0] (119) at (-18.5, -7) {};
		\node [style=new style 0] (120) at (-19.5, -7) {};
		\node [style=new style 0] (121) at (-18, -3) {};
		\node [style=new style 0] (122) at (-21.5, -7) {};
		\node [style=new style 0] (123) at (-22.5, -7) {};
		\node [style=new style 0] (125) at (-22, -3) {};
		\node [style=new style 0] (126) at (-16, 1) {};
		\node [style=new style 0] (127) at (-14.5, -7) {};
		\node [style=new style 0] (128) at (-8.5, -7) {};
		\node [style=new style 0] (129) at (-10.5, -7) {};
		\node [style=new style 0] (130) at (-20.5, -7) {};
		\node [style=new style 0] (131) at (-16.5, -7) {};
		\node [style=new style 0] (132) at (-12.5, -7) {};
		\node [style=new style 0] (133) at (6, -3) {};
		\node [style=new style 0] (134) at (4.5, -7) {};
		\node [style=new style 0] (135) at (6.5, -7) {};
		\node [style=new style 0] (136) at (2.5, -7) {};
		\node [style=new style 0] (137) at (0.5, -7) {};
		\node [style=new style 0] (138) at (2, -3) {};
		\node [style=new style 0] (139) at (-1.5, -7) {};
		\node [style=new style 0] (140) at (-2.5, -7) {};
		\node [style=new style 0] (141) at (-3.5, -7) {};
		\node [style=new style 0] (142) at (-2, -3) {};
		\node [style=new style 0] (143) at (-5.5, -7) {};
		\node [style=new style 0] (144) at (-6.5, -7) {};
		\node [style=new style 0] (145) at (-7.5, -7) {};
		\node [style=new style 0] (146) at (-6, -3) {};
		\node [style=new style 0] (147) at (0, 1) {};
		\node [style=new style 0] (148) at (1.5, -7) {};
		\node [style=new style 0] (149) at (7.5, -7) {};
		\node [style=new style 0] (150) at (5.5, -7) {};
		\node [style=new style 0] (151) at (-4.5, -7) {};
		\node [style=new style 0] (152) at (-0.5, -7) {};
		\node [style=new style 0] (153) at (3.5, -7) {};
		\node [style=new style 0] (154) at (22, -3) {};
		\node [style=new style 0] (155) at (20.5, -7) {};
		\node [style=new style 0] (156) at (22.5, -7) {};
		\node [style=new style 0] (157) at (18.5, -7) {};
		\node [style=new style 0] (158) at (16.5, -7) {};
		\node [style=new style 0] (159) at (18, -3) {};
		\node [style=new style 0] (160) at (14.5, -7) {};
		\node [style=new style 0] (161) at (13.5, -7) {};
		\node [style=new style 0] (162) at (12.5, -7) {};
		\node [style=new style 0] (163) at (14, -3) {};
		\node [style=new style 0] (164) at (10.5, -7) {};
		\node [style=new style 0] (165) at (9.5, -7) {};
		\node [style=new style 0] (166) at (8.5, -7) {};
		\node [style=new style 0] (167) at (10, -3) {};
		\node [style=new style 0] (168) at (16, 1) {};
		\node [style=new style 0] (169) at (17.5, -7) {};
		\node [style=new style 0] (170) at (23.5, -7) {};
		\node [style=new style 0] (171) at (21.5, -7) {};
		\node [style=new style 0] (172) at (11.5, -7) {};
		\node [style=new style 0] (173) at (15.5, -7) {};
		\node [style=new style 0] (174) at (19.5, -7) {};
		\node [style=new style 0] (175) at (-23.5, -7) {};
		\node [style=Black] (176) at (0, 5) {};
		\node [style=Black] (177) at (0, 9) {};
		\node [style=none] (178) at (0, 12) {};
		\node [style=Black] (179) at (22, 5) {};
		\node [style=Black] (180) at (-22, 5) {};
		\node [style=new style 0] (181) at (20.5, 1) {};
		\node [style=new style 0] (182) at (22, 1) {};
		\node [style=new style 0] (184) at (23.5, 1) {};
		\node [style=new style 0] (185) at (-23.5, 1) {};
		\node [style=new style 0] (186) at (-22, 1) {};
		\node [style=new style 0] (188) at (-20.5, 1) {};
		\node [style=none] (189) at (-25, 5) {};
		\node [style=none] (190) at (-25, 1) {};
		\node [style=none] (191) at (-25, -3) {};
		\node [style=none] (192) at (-25, -7) {};
		\node [style=none] (193) at (-25, -11) {};
		\node [style=none] (194) at (-25, 9) {};
		\node [style=none] (195) at (25, -7) {};
		\node [style=none] (196) at (25, -11) {};
		\node [style=none] (197) at (25, -3) {};
		\node [style=none] (198) at (25, 1) {};
		\node [style=none] (199) at (25, 5) {};
		\node [style=none] (200) at (25, 9) {};
		\node [style=none] (201) at (18, -8) {\Huge$\cdots$};
		\node [style=none] (202) at (22, -8) {\Huge$\cdots$};
	\end{pgfonlayer}
	\begin{pgfonlayer}{edgelayer}
		\draw (112) to (113);
		\draw (112) to (114);
		\draw (126) to (125);
		\draw (125) to (123);
		\draw (125) to (122);
		\draw (126) to (121);
		\draw (126) to (117);
		\draw (121) to (120);
		\draw (121) to (119);
		\draw (121) to (118);
		\draw (117) to (115);
		\draw (117) to (127);
		\draw (117) to (116);
		\draw (112) to (129);
		\draw (125) to (130);
		\draw (131) to (121);
		\draw (132) to (117);
		\draw (128) to (112);
		\draw (112) to (126);
		\draw (133) to (134);
		\draw (133) to (135);
		\draw (147) to (146);
		\draw (146) to (145);
		\draw (146) to (144);
		\draw (146) to (143);
		\draw (147) to (142);
		\draw (147) to (138);
		\draw (142) to (141);
		\draw (142) to (140);
		\draw (142) to (139);
		\draw (138) to (136);
		\draw (138) to (148);
		\draw (138) to (137);
		\draw (133) to (150);
		\draw (146) to (151);
		\draw (152) to (142);
		\draw (153) to (138);
		\draw (149) to (133);
		\draw (133) to (147);
		\draw (154) to (155);
		\draw (154) to (156);
		\draw (168) to (167);
		\draw (167) to (166);
		\draw (167) to (165);
		\draw (167) to (164);
		\draw (168) to (163);
		\draw (168) to (159);
		\draw (163) to (162);
		\draw (163) to (161);
		\draw (163) to (160);
		\draw (159) to (157);
		\draw (159) to (169);
		\draw (159) to (158);
		\draw (154) to (171);
		\draw (167) to (172);
		\draw (173) to (163);
		\draw (174) to (159);
		\draw (170) to (154);
		\draw (154) to (168);
		\draw (175) to (125);
		\draw (126) to (176);
		\draw (147) to (176);
		\draw (176) to (168);
		\draw (176) to (177);
		\draw [style=new edge style 0] (177) to (178.center);
		\draw (180) to (177);
		\draw (177) to (179);
		\draw (185) to (180);
		\draw (186) to (180);
		\draw (188) to (180);
		\draw (181) to (179);
		\draw (182) to (179);
		\draw (184) to (179);
		\draw [style=new edge style 0] (200.center) to (194.center);
		\draw [style=new edge style 0] (189.center) to (199.center);
		\draw [style=new edge style 0] (198.center) to (190.center);
		\draw [style=new edge style 0] (191.center) to (197.center);
		\draw [style=new edge style 0] (195.center) to (192.center);
		\draw [style=new edge style 0] (193.center) to (196.center);
	\end{pgfonlayer}
\end{tikzpicture}
  \caption{An example of $\mathfrak{S}_{a,b}$ when $a=3$, $b=4$}\label{Sab}
\end{figure}
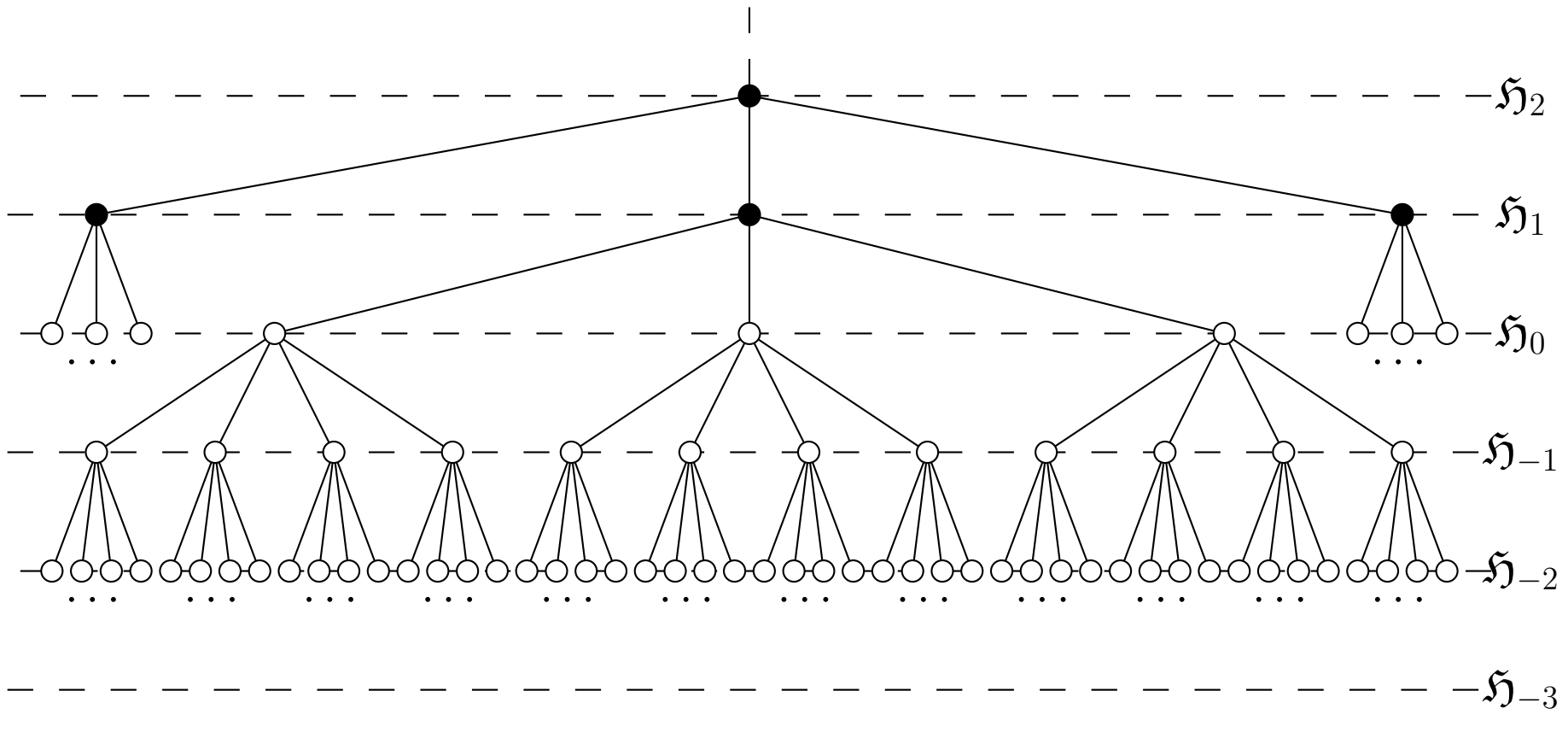

{\bf Example 1:} Motivated by the proof of \cite[Proposition 3.2]{LMSV}, we consider the tree $\mathfrak{S}_{a,b}$ (see Figure \ref{Sab}). For parameters $b \geq a \ge 2$, this tree is defined such that $\nu(x)=b+1$ for all $x \in \mathfrak{H}_j$ with $j \le 0$, and $\nu(x)=a+1$ for all other vertices. For this tree, we have $\mu(T(p^{j}(x),r)) \leq \mu(T(x,r))$ for every $x \in \mathfrak{S}_{a,b}$ and all $j, r \geq 0$.

The construction of $\mathfrak{S}_{a,b}$ leads us to ask: Is the inequality chain \eqref{C1} universally true across the entire class $\Upsilon_{a,b}$? The answer is negative. Indeed, whether \eqref{C1} holds is not determined by the global parameters $(a,b)$ alone but by the intrinsic geometric configuration of the tree. To substantiate this claim, we provide an explicit counterexample in what follows.

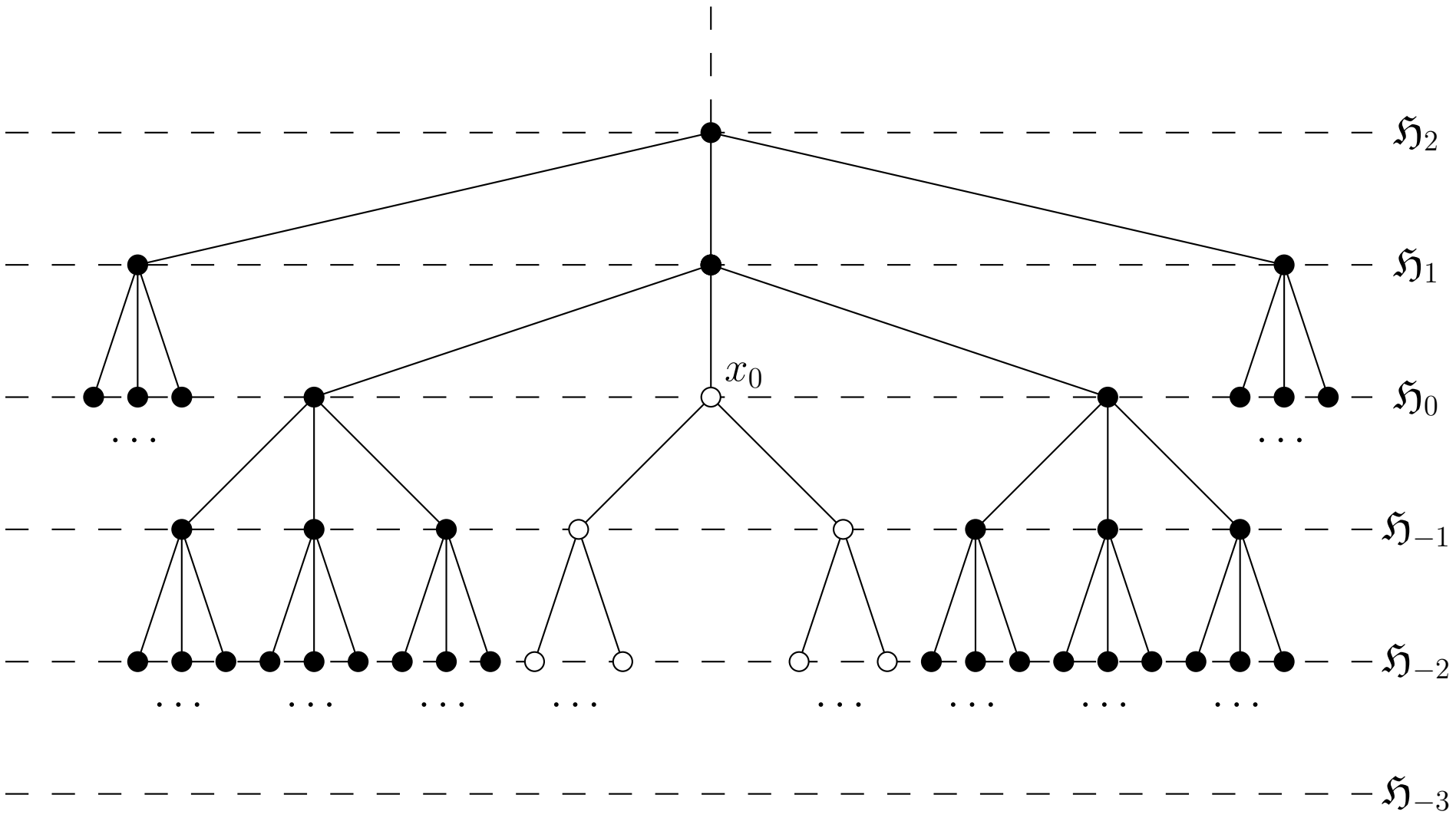
\begin{figure}[!htbp]
  \centering
\tikzstyle{every node}=[scale=0.5]
\begin{tikzpicture}[scale=0.38]
	\begin{pgfonlayer}{nodelayer}
		\node [style=new style 0] (0) at (-3, 0) {};
		\node [style=new style 0] (1) at (0, 3) {};
		\node [style=new style 0] (2) at (3, 0) {};
		\node [style=new style 0] (3) at (-4, -3) {};
		\node [style=new style 0] (4) at (-2, -3) {};
		\node [style=new style 0] (5) at (2, -3) {};
		\node [style=new style 0] (6) at (4, -3) {};
		\node [style=Black] (8) at (0, 6) {};
		\node [style=Black] (9) at (0, 6) {};
		\node [style=Black] (23) at (0, 6) {};
		\node [style=Black] (24) at (9, 3) {};
		\node [style=Black] (25) at (9, 0) {};
		\node [style=Black] (26) at (6, 0) {};
		\node [style=Black] (27) at (6, -3) {};
		\node [style=Black] (28) at (7, -3) {};
		\node [style=Black] (29) at (5, -3) {};
		\node [style=Black] (30) at (8, -3) {};
		\node [style=Black] (31) at (9, -3) {};
		\node [style=Black] (32) at (10, -3) {};
		\node [style=Black] (33) at (11, -3) {};
		\node [style=Black] (34) at (12, -3) {};
		\node [style=Black] (35) at (12, 0) {};
		\node [style=Black] (36) at (13, -3) {};
		\node [style=Black] (37) at (0, 9) {};
		\node [style=Black] (38) at (13, 6) {};
		\node [style=Black] (39) at (-5, -3) {};
		\node [style=Black] (40) at (-7, -3) {};
		\node [style=Black] (41) at (-6, -3) {};
		\node [style=Black] (42) at (-6, 0) {};
		\node [style=Black] (43) at (-8, -3) {};
		\node [style=Black] (44) at (-9, -3) {};
		\node [style=Black] (45) at (-10, -3) {};
		\node [style=Black] (46) at (-9, 0) {};
		\node [style=Black] (47) at (-11, -3) {};
		\node [style=Black] (48) at (-12, -3) {};
		\node [style=Black] (49) at (-13, -3) {};
		\node [style=Black] (50) at (-12, 0) {};
		\node [style=Black] (51) at (-9, 3) {};
		\node [style=Black] (52) at (-13, 6) {};
		\node [style=none] (53) at (0, 12) {};
		\node [style=none] (60) at (-16, 9) {};
		\node [style=none] (61) at (15, 9) {};
		\node [style=none] (62) at (-16, 6) {};
		\node [style=none] (63) at (15, 6) {};
		\node [style=none] (64) at (-16, 3) {};
		\node [style=none] (65) at (15, 3) {};
		\node [style=none] (66) at (-16, 0) {};
		\node [style=none] (67) at (15, 0) {};
		\node [style=none] (68) at (-16, -3) {};
		\node [style=none] (69) at (15, -3) {};
		\node [style=Black] (70) at (-14, 3) {};
		\node [style=Black] (71) at (-13, 3) {};
		\node [style=Black] (72) at (-12, 3) {};
		\node [style=Black] (73) at (12, 3) {};
		\node [style=Black] (74) at (13, 3) {};
		\node [style=Black] (75) at (14, 3) {};
		\node [style=none] (76) at (-16, -6) {};
		\node [style=none] (77) at (15, -6) {};
		\node [style=none] (78) at (0.75, 3.5) {\Huge$x_0$};
		\node [style=none] (79) at (16, 3) {};
		\node [style=none] (80) at (16, 3) {\Huge$\mathfrak{H}_0$};
		\node [style=none] (81) at (16, 6) {\Huge$\mathfrak{H}_1$};
		\node [style=none] (82) at (16, 9) {\Huge$\mathfrak{H}_2$};
		\node [style=none] (83) at (16, 0) {\Huge$\mathfrak{H}_{-1}$};
		\node [style=none] (84) at (16, -3) {\Huge$\mathfrak{H}_{-2}$};
		\node [style=none] (85) at (16, -6) {\Huge$\mathfrak{H}_{-3}$};
		\node [style=none] (87) at (-13, 2) {\Huge$\cdots$};
		\node [style=none] (88) at (13, 2) {\Huge$\cdots$};
		\node [style=none] (89) at (-12, -4) {\Huge$\cdots$};
		\node [style=none] (90) at (-9, -4) {\Huge$\cdots$};
		\node [style=none] (91) at (-6, -4) {\Huge$\cdots$};
		\node [style=none] (92) at (-3, -4) {\Huge$\cdots$};
		\node [style=none] (93) at (3, -4) {\Huge$\cdots$};
		\node [style=none] (94) at (6, -4) {\Huge$\cdots$};
		\node [style=none] (95) at (9, -4) {\Huge$\cdots$};
		\node [style=none] (96) at (12, -4) {\Huge$\cdots$};
	\end{pgfonlayer}
	\begin{pgfonlayer}{edgelayer}
		\draw (0) to (1);
		\draw (1) to (2);
		\draw (9) to (1);
		\draw (0) to (3);
		\draw (0) to (4);
		\draw (2) to (5);
		\draw (2) to (6);
		\draw (51) to (50);
		\draw (50) to (49);
		\draw (50) to (48);
		\draw (50) to (47);
		\draw (51) to (46);
		\draw (51) to (42);
		\draw (46) to (45);
		\draw (46) to (44);
		\draw (46) to (43);
		\draw (42) to (40);
		\draw (42) to (41);
		\draw (42) to (39);
		\draw (51) to (23);
		\draw (37) to (23);
		\draw (37) to (38);
		\draw (37) to (52);
		\draw [style=new edge style 0] (60.center) to (61.center);
		\draw [style=new edge style 0] (62.center) to (63.center);
		\draw [style=new edge style 0] (68.center) to (69.center);
		\draw [style=new edge style 0] (66.center) to (67.center);
		\draw [style=new edge style 0] (64.center) to (65.center);
		\draw [style=new edge style 0] (76.center) to (77.center);
		\draw [style=new edge style 0] (37) to (53.center);
		\draw (29) to (26);
		\draw (27) to (26);
		\draw (28) to (26);
		\draw (30) to (25);
		\draw (31) to (25);
		\draw (32) to (25);
		\draw (33) to (35);
		\draw (34) to (35);
		\draw (36) to (35);
		\draw (35) to (24);
		\draw (25) to (24);
		\draw (26) to (24);
		\draw (24) to (23);
		\draw (70) to (52);
		\draw (71) to (52);
		\draw (72) to (52);
		\draw (75) to (38);
		\draw (74) to (38);
		\draw (73) to (38);
	\end{pgfonlayer}
\end{tikzpicture}
\caption{An example of $\mathfrak{F}_{a,b}$ when $a=2$, $b=3$}\label{Fab}
\end{figure}

{\bf Example 2:} Let $b > a \geq 2$ be integers. The tree $\mathfrak{F}_{a,b}$ (see Figure \ref{Fab}) is defined as an element of $\Upsilon_{a,b}$ for which there exists a distinguished vertex $x_0$ such that $\nu(x_0)=a+1$, $\nu(x)=a+1$ holds for every $x \in \bigcup_{k=0}^{\infty} s^{k}(x_0)$, and $\nu(x)=b+1$ holds for every $x \in \mathfrak{F}_{a,b} \setminus \bigcup_{k=0}^{\infty} s^{k}(x_0)$.

\begin{prop}\label{p3}
Let $\mathfrak{F}_{a,b}$ (with $b > a \geq 2$) be the locally finite tree defined previously. Let $\gamma>0$, $n\in\mathbb{Z}^+$.

(i) Define the function
\begin{align*}
f_{n,2}(y):=\chi_{ T(x_0,n)}(y),\quad y\in\mathfrak{F}_{a,b}.
\end{align*}
Then 
\begin{align*}
\lim_{n\to+\infty}
\frac{\mathcal{N}^{\gamma}f_{n,2}(x_0)}{\mathcal{M}^{\gamma}f_{n,2}(x_0)}
=+\infty.
\end{align*}

(ii) If $a<b<a^2$ and $\gamma\geq\frac{1}{2\log_b a-1}$, then
\begin{align*}
\mathcal{N}^{\gamma}f(x)
\lesssim\widetilde{\mathcal{M}}^{\gamma}f(x)
\end{align*}
for every $x\in\mathfrak{F}_{a,b}$.
\end{prop}

\begin{proof}
(i) A straightforward calculation shows that 
\begin{align*}
\mu(B(x_0,n))=\frac{a^{n+1}-1}{a-1}+\frac{b^{n}-1}{b-1},\quad n\in\mathbb{Z}^+,
\end{align*}
and
\begin{align*}
\mu( T(x_0,n))=\frac{a^{n+1}-1}{a-1},\quad n\in\mathbb{Z}^+.
\end{align*}
Thus we have that
\begin{align*}
\frac{\mathcal{N}^{\gamma}f_{n,2}(x_0)}{\mathcal{M}^{\gamma}f_{n,2}(x_0)}
=\frac{\frac{\mu( T(x_0,n))}{\mu(T(x_0,n))^\gamma}}{\frac{\mu( T(x_0,n))}{\mu(B(x_0,n))^\gamma}}
>\Big(\frac{a-1}{b-1}\frac{b^{n}-1}{a^{n+1}-1}\Big)^{\gamma}\to+\infty,\; \text{as}\; n\to+\infty. 
\end{align*} 
It follows that $\mathcal{N}^{\gamma}$ and $\widetilde{\mathcal{N}}^{\gamma}$ are not uniformly dominated by $\mathcal{M}^{\gamma}$ for all $f$ on $\mathfrak{F}_{a,b}$.

(ii) Note that $a^r\lesssim\mu(B(x,r))\lesssim b^r$ for any $x\in\mathfrak{F}_{a,b}$ and $r\in\mathbb{N}$. Then for every $x\in\mathfrak{F}_{a,b}$, we have that 
\begin{align*}
\frac{1}{\mu(T(x,r))^\gamma}\int_{T(x,r)}\vert f\vert d\mu
&\leq\frac{1}{\mu(T(x,r))^\gamma}\int_{\bigcup\limits_{y\in s^{\lfloor r/2\rfloor}(x)}B(y,\lfloor r/2\rfloor+1)}\vert f\vert d\mu\\
&\leq\sum_{y\in s^{\lfloor r/2\rfloor}(x)}\frac{1}{\mu(T(x,r))^\gamma}\int_{B(y,\lfloor r/2\rfloor+1)}\vert f\vert d\mu\\
&\leq\sum_{y\in s^{\lfloor r/2\rfloor}(x)}\Big(\frac{\mu(B(y,\lfloor r/2\rfloor+1))}{\mu(T(x,r))}\Big)^\gamma
\widetilde{\mathcal{M}}^{\gamma}f(x)\\
&\leq\frac{\mu(s^{\lfloor r/2\rfloor}(x))\max\limits_{y\in s^{\lfloor r/2\rfloor}(x)}\mu(B(y,\lfloor r/2\rfloor+1))^\gamma}{\mu(T(x,r))^\gamma}\widetilde{\mathcal{M}}^{\gamma}f(x)\\
&\lesssim\Big(\frac{b^{1+\gamma}}{a^{2\gamma}}\Big)^{\lfloor r/2\rfloor}\widetilde{\mathcal{M}}^{\gamma}f(x).
\end{align*}

Since $\gamma\geq\frac{1}{2\log_b a-1}$, it follows that $b^{1+\gamma}\leq a^{2\gamma}$. Taking the supremum over $r\in\mathbb{N}$, we have
\begin{align*}
\mathcal{N}^{\gamma}f(x)
\lesssim\widetilde{\mathcal{M}}^{\gamma}f(x).
\end{align*}
This completes the proof.
\end{proof}

Two natural questions remain open for the class $\Upsilon_{a,b}$ and any $\gamma>0$: (i) whether $\mathcal{N}^{\gamma}$ and $\widetilde{\mathcal{N}}^{\gamma}$ are pointwise dominated by $\widetilde{\mathcal{M}}^{\gamma}$; (ii) whether $\mathcal{M}^{\gamma}$ is pointwise dominated by $\widetilde{\mathcal{N}}^{\gamma}$. However, in the setting of unbounded geometry, we demonstrate by construction that there exists a class of trees on which both $\widetilde{\mathcal{N}}^{\gamma}\lesssim\widetilde{\mathcal{M}}^{\gamma}$ and $\mathcal{M}^{\gamma}\lesssim\widetilde{\mathcal{N}}^{\gamma}$ fail.

\begin{figure}[!htbp]
  \centering
\tikzstyle{every node}=[scale=0.5]
\begin{tikzpicture}[scale=0.37]
	\begin{pgfonlayer}{nodelayer}
		\node [style=new style 0] (0) at (-7.5, 0) {};
		\node [style=new style 0] (1) at (-6.5, 3) {};
		\node [style=new style 0] (2) at (-5.5, 0) {};
		\node [style=Black] (23) at (-6.5, 6) {};
		\node [style=new style 0] (24) at (-2.5, 3) {};
		\node [style=new style 0] (26) at (-3.5, 0) {};
		\node [style=new style 0] (35) at (-1.5, 0) {};
		\node [style=Black] (37) at (0, 9) {};
		\node [style=Black] (38) at (13, 6) {};
		\node [style=new style 0] (42) at (-9.5, 0) {};
		\node [style=new style 0] (50) at (-11.5, 0) {};
		\node [style=new style 0] (51) at (-10.5, 3) {};
		\node [style=Black] (52) at (-13, 6) {};
		\node [style=Black] (53) at (0, 12) {};
		\node [style=none] (60) at (-16, 9) {};
		\node [style=none] (61) at (16, 9) {};
		\node [style=none] (62) at (-16, 6) {};
		\node [style=none] (63) at (16, 6) {};
		\node [style=none] (64) at (-16, 3) {};
		\node [style=none] (65) at (16, 3) {};
		\node [style=none] (66) at (-16, 0) {};
		\node [style=none] (67) at (16, 0) {};
		\node [style=none] (68) at (-16, -3) {};
		\node [style=none] (69) at (16, -3) {};
		\node [style=new style 0] (70) at (-14, 3) {};
		\node [style=new style 0] (71) at (-13, 3) {};
		\node [style=new style 0] (72) at (-12, 3) {};
		\node [style=new style 0] (73) at (12, 3) {};
		\node [style=new style 0] (74) at (13, 3) {};
		\node [style=new style 0] (75) at (14, 3) {};
		\node [style=none] (79) at (17, 3) {};
		\node [style=none] (80) at (17, 3) {\Huge$\mathfrak{H}_0$};
		\node [style=none] (81) at (17, 6) {\Huge$\mathfrak{H}_1$};
		\node [style=none] (82) at (17, 9) {\Huge$\mathfrak{H}_2$};
		\node [style=none] (83) at (17, 0) {\Huge$\mathfrak{H}_{-1}$};
		\node [style=none] (84) at (17, -3) {\Huge$\mathfrak{H}_{-2}$};
		\node [style=none] (87) at (-13, 2) {\Huge$\cdots$};
		\node [style=none] (88) at (13, 2) {\Huge$\cdots$};
		\node [style=none] (89) at (-11.5, -1) {\Huge$\cdots$};
		\node [style=new style 0] (97) at (5.5, 0) {};
		\node [style=new style 0] (98) at (6.5, 3) {};
		\node [style=new style 0] (99) at (7.5, 0) {};
		\node [style=Black] (105) at (6.5, 6) {};
		\node [style=new style 0] (107) at (10.5, 3) {};
		\node [style=new style 0] (108) at (9.5, 0) {};
		\node [style=new style 0] (112) at (11.5, 0) {};
		\node [style=new style 0] (116) at (3.5, 0) {};
		\node [style=new style 0] (119) at (1.5, 0) {};
		\node [style=new style 0] (120) at (2.5, 3) {};
		\node [style=none] (121) at (-16, 12) {};
		\node [style=none] (122) at (16, 12) {};
		\node [style=none] (123) at (17, 12) {\Huge$\mathfrak{H}_{3}$};
		\node [style=Black] (124) at (-14, 9) {};
		\node [style=Black] (125) at (-7, 9) {};
		\node [style=Black] (126) at (7, 9) {};
		\node [style=Black] (127) at (14, 9) {};
		\node [style=none] (128) at (0, 15) {};
		\node [style=none] (129) at (-9.5, -1) {\Huge$\cdots$};
		\node [style=none] (130) at (-7.5, -1) {\Huge$\cdots$};
		\node [style=none] (131) at (-5.5, -1) {\Huge$\cdots$};
		\node [style=none] (132) at (-3.5, -1) {\Huge$\cdots$};
		\node [style=none] (133) at (-1.5, -1) {\Huge$\cdots$};
		\node [style=none] (134) at (1.5, -1) {\Huge$\cdots$};
		\node [style=none] (135) at (3.5, -1) {\Huge$\cdots$};
		\node [style=none] (136) at (5.5, -1) {\Huge$\cdots$};
		\node [style=none] (137) at (7.5, -1) {\Huge$\cdots$};
		\node [style=none] (138) at (9.5, -1) {\Huge$\cdots$};
		\node [style=none] (139) at (11.5, -1) {\Huge$\cdots$};
		\node [style=none] (140) at (-14, 8) {\Huge$\cdots$};
		\node [style=none] (141) at (14, 8) {\Huge$\cdots$};
		\node [style=none] (142) at (7, 8) {\Huge$\cdots$};
		\node [style=none] (143) at (-7, 8) {\Huge$\cdots$};
	\end{pgfonlayer}
	\begin{pgfonlayer}{edgelayer}
		\draw (0) to (1);
		\draw (1) to (2);
		\draw (51) to (50);
		\draw (51) to (42);
		\draw (51) to (23);
		\draw (37) to (23);
		\draw (37) to (38);
		\draw (37) to (52);
		\draw [style=new edge style 0] (60.center) to (61.center);
		\draw [style=new edge style 0] (62.center) to (63.center);
		\draw [style=new edge style 0] (68.center) to (69.center);
		\draw [style=new edge style 0] (66.center) to (67.center);
		\draw [style=new edge style 0] (64.center) to (65.center);
		\draw (37) to (53);
		\draw (35) to (24);
		\draw (26) to (24);
		\draw (24) to (23);
		\draw (70) to (52);
		\draw (71) to (52);
		\draw (72) to (52);
		\draw (75) to (38);
		\draw (74) to (38);
		\draw (73) to (38);
		\draw (97) to (98);
		\draw (98) to (99);
		\draw (105) to (98);
		\draw (120) to (119);
		\draw (120) to (116);
		\draw (112) to (107);
		\draw (108) to (107);
		\draw [style=new edge style 0] (121.center) to (122.center);
		\draw (124) to (53);
		\draw (53) to (125);
		\draw (53) to (126);
		\draw (53) to (127);
		\draw [style=new edge style 0] (128.center) to (53);
		\draw (120) to (105);
		\draw (105) to (107);
		\draw (23) to (1);
		\draw (37) to (105);
	\end{pgfonlayer}
\end{tikzpicture}
\caption{An example of $\mathfrak{P}_2$}\label{P2}
\end{figure}
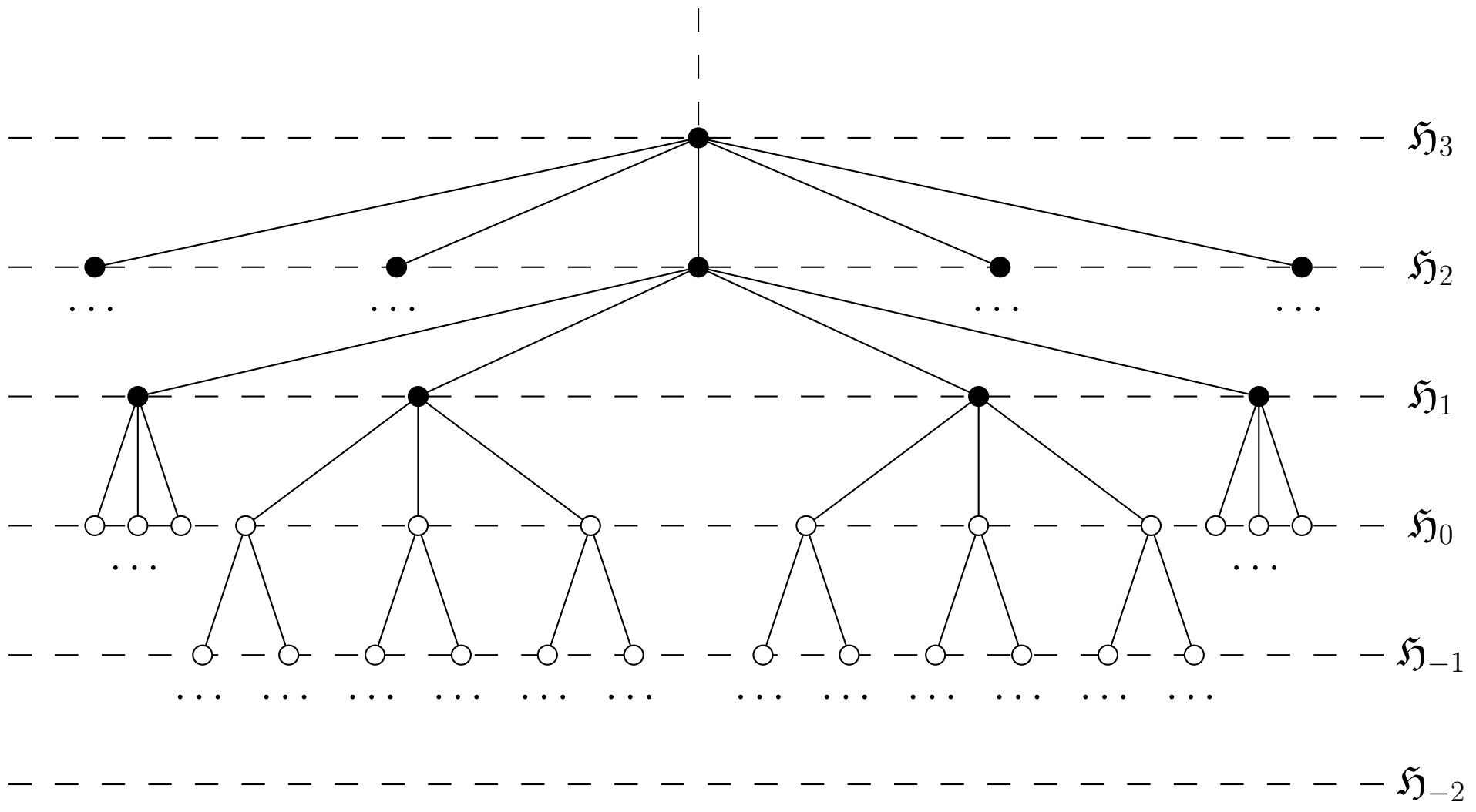

{\bf Example 3:} Let $a\geq2$ and $a\in\mathbb{N}$. The tree $\mathfrak{P}_a$ (see Figure \ref{P2}) has a degree that grows linearly with the height from a fixed base level. Formally, its degree function is defined as:
\begin{equation*}
\nu(x)=
\begin{cases}
a+1 & {h(x)\leq0,}\\
a+j+1 & {h(x)=j\geq1.}
\end{cases}
\end{equation*}

\begin{prop}\label{p5}
Let  $\gamma>0$, $n\in\mathbb{Z}^+$, $x_0$ be the reference vertex of $\mathfrak{P}_a$, $\omega=\{x_0, p(x_0), p^2(x_0), \ldots \}$ be the geodesic ray of $\mathfrak{P}_a$ and $z_n\in s^n(p^n(x_0))$ such that $h(x_0)=0$, $z_n\wedge x_0=p^{n}(x_0)$ and $d(z_n,x_0)=2n$. Then 
\begin{align}\label{P2D}
\lim_{n\to+\infty}
\frac{\widetilde{\mathcal{N}}^{\gamma}\delta_{z_n}(x_0)}
{\widetilde{\mathcal{M}}^{\gamma}\delta_{z_n}(x_0)}
=\lim_{n\to+\infty}\frac{\mathcal{M}^{\gamma}\delta_{p^n(x_0)}(x_0)}
{\widetilde{\mathcal{N}}^{\gamma}\delta_{p^n(x_0)}(x_0)}=+\infty.
\end{align}
\end{prop}
\begin{proof}
We first prove that the first term in \eqref{P2D} diverges. It follows from 
\begin{align*}
\widetilde{\mathcal{N}}^{\gamma}\delta_{z_n}(x_0)
=\sup_{k\in\mathbb{N}}\sup_{r\geq0}
\frac{1}{\mu(T(p^{k}(x_0),r))^{\gamma}}\int_{T(p^{k}(x_0),r)}\delta_{z_n}(y) d\mu(y)
\end{align*}
that $\widetilde{\mathcal{N}}^{\gamma}\delta_{z_n}(x_0)>0$ if and only if there exist some $k\in\mathbb{N}$ such that $p^{k}(x_0)\succeq z_n$. It is easy to see that $z_n\in s^k(p^k(x_0))$ if and only if $k\geq n$. Then we obtain that
\begin{align*}
\widetilde{\mathcal{N}}^{\gamma}\delta_{z_n}(x_0)
&=\sup_{k\geq n}\sup_{r\geq n}
\frac{1}{\mu(T(p^{k}(x_0),r))^{\gamma}}\int_{T(p^{k}(x_0),r)}\delta_{z_n}(y) d\mu(y)\\
&={\mu(T(p^{n}(x_0),n))^{-\gamma}}.
\end{align*}

Let $B(y,r_y)$ be the ball containing both $x_0$ and $z_n$ for any $y\in\mathfrak{P}_a$.
For any vertex $v\in[x_0,z_n]$, we consider the geodesic $[y,v]$ from $y$ to $v$. We may write it as 
\[y=u_0\sim u_1\sim\cdots\sim u_k=v,\; d(y,u_i)=i\; (0\leq i\leq k).\]
Since $v\in[x_0,z_n]$, the set $\{u_i\in[y,v]:u_i\in[x_0,z_n]\}$ is nonempty and finite.
Hence there exists a minimal index $j$ such that $u_j\in[x_0,z_n]$.
Without loss of generality, we can assume that $u_j\in[x_0,x_0\wedge z_n]$.
The vertex $u_j$ is precisely the first vertex where the paths $[y,v]$ and $[x_0,z_n]$ intersect. Then $d(y,v)=d(y,u_j)+d(u_j,v)$. It follows from $d(u_j,v)\leq\max\{d(u_j,x_0),d(u_j,z_n)\}$ that
\begin{align*}
d(y,v)&\leq\max\{d(y,u_j)+d(u_j,x_0),d(y,u_j)+d(u_j,z_n)\}\\
&=\max\{d(y,x_0),d(y,z_n)\}\\
&\leq r_y.
\end{align*}
Hence $[x_0,z_n]\subset B(y,r_y)$ and $x_0\wedge z_n\in B(y,r_y)$. Now, for any vertex $z\in B(x_0\wedge z_n,n)$, we have that
\begin{align*}
d(z,y)&\leq d(z,x_0\wedge z_n)+d(x_0\wedge z_n,u_j)+d(u_j,y)\\
&\leq n+d(x_0\wedge z_n,u_j)+d(u_j,y)\\
&=d(x_0\wedge z_n, z_n)+d(x_0\wedge z_n,u_j)+d(u_j,y)\\
&=d(z_n,u_j)+d(u_j,y)\\
&=d(y,z_n)\leq r_y.
\end{align*} 
Therefore, any ball $B(y, r_y)$ that contains both $x_0$ and $z_n$ must also contain $B(x_0 \wedge z_n, n)$. Then we have 
\[\widetilde{\mathcal{M}}^{\gamma}\delta_{z_n}(x_0)={\mu(B(x_0\wedge z_n,n))^{-\gamma}}.\]
By Lemma \ref{BASE2}, as $n\to+\infty$, we obtain that
\begin{align*}
\frac{\widetilde{\mathcal{N}}^{\gamma}\delta_{z_n}(x_0)}
{\widetilde{\mathcal{M}}^{\gamma}\delta_{z_n}(x_0)}
&=\Big(\frac{\mu(B(x_0\wedge z_n,n))}{\mu(T(x_0\wedge z_n,n))}\Big)^\gamma\\
&\gtrsim\Big(\frac{\mu(s^{n-1}(p(x_0\wedge z_n)))}{\mu(\beta(T(x_0\wedge z_n,n)))}\Big)^\gamma\\
&\gtrsim\Big(\frac{(a+n+1)!}{(a+n)!}\Big)^\gamma\to+\infty.
\end{align*}

Now we estimate the second term of \eqref{P2D}. Without loss of generality, we let $n$ be an even integer, and write $n = 2k$ with $k \in \mathbb{Z}^+$. Then we have that 
\begin{align*}
\frac{\mathcal{M}^{\gamma}\delta_{p^n(x_0)}(x_0)}
{\widetilde{\mathcal{N}}^{\gamma}\delta_{p^n(x_0)}(x_0)}
&=
\frac{\displaystyle\frac{1}{\mu(B(x_0,n))^\gamma}
\int_{B(x_0,n)}\delta_{p^n(x_0)}(y)d\mu(y)}
{\displaystyle\frac{1}{\mu(T(p^n(x_0),n))^\gamma}\int_{ T(p^n(x_0),n)}\delta_{p^n(x_0)}(y)d\mu(y)}
\\
&=\Big(\frac{\mu(T(p^n(x_0),n))}{\mu(B(x_0,n))}\Big)^\gamma.
\end{align*}
From Lemma \ref{BASE2} (ii) and \eqref{MUBUT}, it follows that 
\[
\mu(T(p^n(x_0),n))\approx\mu(\beta(T(p^n(x_0),n)))\approx(a+n)!
\]
and \\[-1.5em]
\begin{align*}
\mu(B(x_0,n))
&\approx\sum_{j=0}^{n}\mu(\beta(T(p^j(x_0),n-j)))\\
&=a^n+\sum_{j=1}^{k-1}\frac{(a+j)!}{a!}a^{n-j-1}+\sum_{j=k}^{n}\binom{a+j}{n-j}(n-j)!\\
&\leq a^n+(k-1)\frac{(a+k)!}{a!}a^{n-2}+\sum_{j=0}^{k}(a+n)^j\\
&=:a^n+S_1+S_2.
\end{align*}

We claim that $a^n+S_1+S_2=o((a+n)!)$. 

It is trivial to obtain that 
\[
\lim\limits_{n\to+\infty}\frac{a^n}{(a+n)!}=0.
\] 

For $S_1$, we have that
\begin{align*}
\frac{S_1}{(a+n)!}=\frac{(k-1)a^{n-2}}{a!}\frac{1}{\prod_{i=1}^{n-k}(a+k+i)}.
\end{align*}
Take $n \geq 4a^2$ sufficiently large and write $n = 2k$. As $n \to +\infty$, we have 
\begin{align*}
\frac{S_1}{(a+n)!}=\frac{(k-1)a^{n-2}}{a!\prod_{i=1}^{n-k}(a+k+i)}
\leq\frac{(k-1)a^{n-2}}{a!}\Big(\frac{2}{n}\Big)^{\frac{n}{2}}
\lesssim n\Big(\frac{2a^2}{n}\Big)^{\frac{n}{2}}
\to 0.
\end{align*} 

For $S_2$, we have that 
\begin{align*}
\frac{S_2}{(a+n)!}=\frac{(a+2k)^{k+1}-1}{(a+2k-1)\cdot(a+2k)!}
\approx\frac{(a+2k)^k}{(a+2k)!}.
\end{align*}
Denote by $u_k:=\frac{(a+2k)^k}{(a+2k)!}$.
It follows that 
\begin{align*}
0<\frac{u_{k+1}}{u_k}
=\Big(1+\frac{2}{a+2k}\Big)^{k}\frac{1}{a+2k+1}
<\left(1+\frac{1}{k}\right)^{k}\frac{1}{2k}
\to 0,
\end{align*}
as $k\to+\infty$. Hence $u_k\to 0$, as $k\to+\infty$. This completes the proof.
\end{proof}

\begin{longtable}[c]{c|c|c|c|c}
$\lesssim$ &
  $\mathcal{M}^\gamma$ &
  $\widetilde{\mathcal{M}}^\gamma$ &
  $\mathcal{N}^\gamma$ &
  $\widetilde{\mathcal{N}}^\gamma$ \\ \hline
\endfirsthead
\endhead
\hline
\endfoot
\endlastfoot
\multirow{2}{*}{$\mathcal{M}^\gamma$} &
  \multirow{2}{*}{-} &
  \multirow{2}{*}{Yes} &
  No &
  Yes \\
 &
   &
   &
  Prop. \ref{phc} (i) &
  Prop. \ref{p1} \\ \hline
\multirow{2}{*}{$\widetilde{\mathcal{M}}^\gamma$} &
  No &
  \multirow{2}{*}{-} &
  No &
  No \\
 &
  Prop. \ref{phc} (ii) &
   &
  Prop. \ref{phc} (i) &
  Prop. \ref{phc} (ii) \\ \hline
\multirow{2}{*}{$\mathcal{N}^\gamma$} &
  Yes &
  Yes &
  \multirow{2}{*}{-} &
  \multirow{2}{*}{Yes} \\
 &
  Prop. \ref{p1} &
  Prop. \ref{p1} &
   &
   \\ \hline
\multirow{2}{*}{$\widetilde{\mathcal{N}}^\gamma$} &
  No &
  Yes &
  No &
  \multirow{2}{*}{-} \\
 &
  Prop. \ref{phc} (iii) &
  Prop. \ref{p1} &
  Prop. \ref{phc} (i) &
   \\ \hline
\caption{Pointwise dominance relations among $\mathcal{N}^{\gamma}$,
$\widetilde{\mathcal{N}}^{\gamma}$, $\mathcal{M}^{\gamma}$ and $\widetilde{\mathcal{M}}^{\gamma}$ on homogeneous trees}
\label{T1}
\end{longtable}


\begin{longtable}[c]{c|c|c|c|c}
$\lesssim$ & $\mathcal{M}^\gamma$                  & $\widetilde{\mathcal{M}}^\gamma$ & $\mathcal{N}^\gamma$                & $\widetilde{\mathcal{N}}^\gamma$     \\ \hline
\endfirsthead
\endhead
\hline
\endfoot
\endlastfoot
\multirow{2}{*}{$\mathcal{M}^\gamma$} &
  \multirow{2}{*}{-} &
  \multirow{2}{*}{Yes} &
  No &
  No \\
 &
   &
   &
  Prop. \ref{phc} (i) &
  Prop. \ref{p5} \\ \hline
\multirow{2}{*}{$\widetilde{\mathcal{M}}^\gamma$} &
  No &
  \multirow{2}{*}{-} &
  No &
  No \\
           & Prop. \ref{phc} (ii)  &                                  & Prop. \ref{phc} (i) & Prop. \ref{phc} (ii) \\ \hline
\multirow{2}{*}{$\mathcal{N}^\gamma$} &
  No &
\makecell[c]{Yes on some $\mathfrak{T}\in\Upsilon_{a,b}$\\ and $\gamma\geq(2\log_b a-1)^{-1}$} &
  \multirow{2}{*}{-} &
  \multirow{2}{*}{Yes} \\
 &
  Prop. \ref{p3} (i) &
  Prop. \ref{p3} (ii) &
   &
   \\ \hline
\multirow{2}{*}{$\widetilde{\mathcal{N}}^\gamma$} &
  No &
  No &
  No &
  \multirow{2}{*}{-} \\
           & Prop. \ref{phc} (iii) & Prop. \ref{p5}  & Prop. \ref{phc} (i) &                                      \\ \hline

\caption{Pointwise dominance relations among $\mathcal{N}^{\gamma}$,
$\widetilde{\mathcal{N}}^{\gamma}$, $\mathcal{M}^{\gamma}$ and $\widetilde{\mathcal{M}}^{\gamma}$ on locally finite trees}
  \label{T2}
\end{longtable}

Consequently, there is no universal set of pointwise inequalities among these four operators that holds for all functions $f$. Whether one operator yields a larger value than another at a given point can vary, depending on both the local geometry of the tree and the specific profile of the function $f$. We summarize the above results in the tables above.

\section{Boundedness and Optimality of Fractional Triangular Maximal Operators} \label{s3}

We establish two main results regarding the operators $\mathcal{N}^\gamma$ and $\widetilde{\mathcal{N}}^\gamma$: their strong-type boundedness from $\ell^p$ to $\ell^q$ on locally finite trees, and the optimality of the exponent ranges.

\begin{thm}\label{thm1}
Let $0<\gamma\leq1$. Suppose that $\mathfrak{T}$ is a locally finite tree with $3\leq a+1\leq\nu<+\infty$. Then 

(i) $\mathcal{N}^\gamma$ is of strong type $(p,q)$ if $1\leq p\leq q\leq+\infty$ and $1\leq p<\frac{1}{1-\gamma}$;

(ii) $\mathcal{N}^{\gamma}$ is unbounded from $\ell^{p}(\mathfrak{T})$ to $\ell^{q}(\mathfrak{T})$ if $\frac{1}{1-\gamma}< p\leq q\leq+\infty$ or $p=q=\frac{1}{1-\gamma}$.
\end{thm}

\begin{proof}
(i) First, we prove that $\mathcal{N}^\gamma$ is of strong type $(1,1)$. Denote by $E:=\{(x,y)\in\mathfrak{T}\times\mathfrak{T}:x\succeq y\}$. For $x\in\mathfrak{T}$ and $r>0$, we observe that if $y\in T(x,r)$, then $ T(x,d(x,y))\subset T(x,r)$. It follows that
\begin{align*}
\frac{1}{\mu( T(x,r))^{\gamma}}\int_{ T(x,r)}\vert f(y)\vert d\mu(y)
&\leq\int_{ T(x,r)}\frac{1}{\mu( T(x,d(x,y)))^\gamma}\vert f(y)\vert d\mu(y)\\
&\leq\int_{\{y\in\mathfrak{T}:x\succeq y\}}\frac{1}{\mu( T(x,d(x,y)))^\gamma}\vert f(y)\vert d\mu(y)\\
&=\int_{\mathfrak{T}}\frac{\chi_{E}(x,y)}{\mu( T(x,d(x,y)))^\gamma}\vert f(y)\vert d\mu(y)\\
&=:\int_{\mathfrak{T}}\tau_{\gamma}(x,y)\vert f(y)\vert d\mu(y).
\end{align*}
Thus, we have that
\begin{align}\label{tau}
\mathcal{N}^\gamma f(x)\leq\int_{\mathfrak{T}}\tau_{\gamma}(x,y)\vert f(y)\vert d\mu(y)
\end{align}
and 
\begin{align*}
\Vert \mathcal{N}^\gamma f\Vert_{1}
\leq \int_{\mathfrak{T}}\int_{\mathfrak{T}}\tau_{\gamma}(x,y)\vert f(y)\vert d\mu(y)d\mu(x)
\leq\Vert f\Vert_{1}\int_{\mathfrak{T}}\sup_{y\in\mathfrak{T}}\tau_{\gamma}(x,y)d\mu(x).
\end{align*}

We claim that $\sup_{y\in\mathfrak{T}}\int_{\mathfrak{T}}\tau_{\gamma}(x,y)d\mu(x)<+\infty$. In fact, given $y\in\mathfrak{T}$, if $x\succeq y$, then exists $k\in\mathbb{N}$ such that $x=p^{k}(y)$. Then we have
\[\sup_{y\in\mathfrak{T}}\int_{\mathfrak{T}}\tau_{\gamma}(x,y)d\mu(x)
=\sup_{y\in\mathfrak{T}}\sum_{k=0}^{+\infty}\frac{1}{\mu( T(p^{k}(y),k))^\gamma}.\]
Since 
\[\mu( T(p^{k}(y),k))^\gamma
=\Big(\sum_{j=0}^{k}\mu(s^{j}(p^{k}(y)))\Big)^{\gamma}
\geq\Big(\sum_{j=0}^{k}a^j\Big)^{\gamma}
\geq a^{k\gamma}\]
for every $y\in\mathfrak{T}$, it follows that 
\[\sup_{y\in\mathfrak{T}}\int_{\mathfrak{T}}\tau_{\gamma}(x,y)d\mu(x)
\leq\sum_{k=0}^{\infty}\frac{1}{a^{k\gamma}}
=\frac{a^\gamma}{a^{\gamma}-1}<+\infty.\]
This proves that 
$\Vert \mathcal{N}^{\gamma}f\Vert_1\leq\frac{a^\gamma}{a^{\gamma}-1}\Vert f\Vert_1$ for all $f\in\ell^1(\mathfrak{T})$. Since the continuous inclusion $\ell^{p_1}(\mathfrak{T})\hookrightarrow\ell^{p_2}(\mathfrak{T})$ holds if $1\leq p_1\leq p_2\leq+\infty$. Thus we obtain that $\mathcal{N}^{\gamma}$ is of strong type $(1,q)$ for every $1\leq q\leq+\infty$.

Second, we prove that $\mathcal{N}^\gamma$ is of strong type $(p,p)$ if $1<p<{1}/{(1-\gamma)}$. If $1<p<{1}/{(1-\gamma)}$, then $p'>{1}/{\gamma}$ where $p'$ is the conjugate index of $p$. Choose $0<\beta<1$ such that  ${1}/{\gamma}<p'\beta<p'$.  By H\"{o}lder's inequality, it follows that
\begin{align*}
\int_{\mathfrak{T}}\tau_{\gamma}(x,y)\vert f(y)\vert d\mu(y)
\!\leq\!\Big(\int_{\mathfrak{T}}\tau_{\gamma}(x,y)^{p(1-\beta)}\vert f(y)\vert^p d\mu(y)\Big)^{\frac{1}{p}}
\Big(\int_{\mathfrak{T}}\tau_{\gamma}(x,y)^{p'\beta}d\mu(y)\Big)^{\frac{1}{p'}}.
\end{align*}

Since $s^k(x)\subset T(x,k)$, it follows that 
\begin{align*}
\int_{\mathfrak{T}}\tau_{\gamma}(x,y)^{p'\beta} d\mu(y)
&=\sum_{k=0}^{+\infty}\frac{\mu(s^k(x))}{\mu( T(x,k))^{p'\beta\gamma}}\\
&\leq\sum_{k=0}^{+\infty}\frac{1}{\mu( T(x,k))^{p'\beta\gamma-1}}\\
&\leq\sum_{k=0}^{+\infty}\frac{1}{a^{k(p'\beta\gamma-1)}}\\
&<+\infty.
\end{align*}
This implies that
\begin{align*}
\int_{\mathfrak{T}}\tau_{\gamma}(x,y)\vert f(y)\vert d\mu(y)
\lesssim\Big(\int_{\mathfrak{T}}\tau_{\gamma}(x,y)^{p(1-\beta)}\vert f(y)\vert^p d\mu(y)\Big)^{{1}/{p}}.
\end{align*}
Thus, by \eqref{tau}, we have that 
\begin{align*}
\Vert \mathcal{N}^\gamma f\Vert_{p}^p
&\leq\int_{\mathfrak{T}}\Big(\int_{\mathfrak{T}}\tau_{\gamma}(x,y)\vert f(y)\vert d\mu(y)\Big)^p d\mu(x)\\
&\lesssim\int_{\mathfrak{T}}\int_{\mathfrak{T}}\tau_{\gamma}(x,y)^{p(1-\beta)}\vert f(y)\vert^p d\mu(y)d\mu(x)\\
&\leq\Vert f\Vert_{p}^{p}
\sup_{y\in\mathfrak{T}}\int_{\mathfrak{T}}\tau_{\gamma}(x,y)^{p(1-\beta)} d\mu(x).
\end{align*}
Here, similar as the above estimate,
\[\sup_{y\in\mathfrak{T}}\int_{\mathfrak{T}}\tau_{\gamma}(x,y)^{p(1-\beta)}d\mu(x)
=\sup_{y\in\mathfrak{T}}
\sum_{k=0}^{+\infty}\frac{1}{\mu( T(p^{k}(y),k))^{p\gamma(1-\beta)}}.\]
Hence we have that
\[\sup_{y\in\mathfrak{T}}
\int_{\mathfrak{T}}\tau_{\gamma}(x,y)^{p(1-\beta)}d\mu(x)<+\infty,\]
which implies that $\mathcal{N}^\gamma$ is of strong type $(p,p)$ if $1<p<{1}/{(1-\gamma)}$.

Using H\"older's inequality, for every $x\in\mathfrak{T}$ and $r\geq0$, we have that
\begin{align}\label{Holder}
\frac{1}{\mu( T(x,r))^{\gamma}}\int_{ T(x,r)}\vert f(y)\vert d\mu
\leq\Vert f\Vert_{\frac{1}{1-\gamma}},
\end{align}
which yields that
\[\Vert \mathcal{N}^{\gamma}f\Vert_{\infty}\leq\Vert f\Vert_{\frac{1}{1-\gamma}}.\]

Therefore, the Marcinkiewicz interpolation theorem implies that $\mathcal{N}^{\gamma}$ is of strong type $(p,q)$ if $1\leq p\leq q\leq+\infty$ and $1\leq p<{1}/{(1-\gamma)}$.

(ii) For $n \in \mathbb{Z}^+$ and a fixed vertex $o \in \mathfrak{T}$, we consider the following cases.

Case 1: For ${1}/{(1-\gamma)} < p \leq q \leq +\infty$, set $g_{n,1}(x) := \chi_{T(o,n)}(x)$. Then we have 
\begin{align*}
\mathcal{N}^{\gamma}g_{n,1}(o)
=\frac{\mu(T(o,n))}{\mu( T(o,n))^{\gamma}}.
\end{align*}

If ${1}/{(1-\gamma)}<p<\infty$, then we have
\begin{align*}
\frac{\Vert \mathcal{N}^{\gamma}g_{n,1}\Vert_q}{\Vert g_{n,1}\Vert_p}
\geq\frac{\mathcal{N}^{\gamma}g_{n,1}(o)}{\mu(T(o,n))^{\frac{1}{p}}}
\geq\mu(T(o,n))^{1-\gamma-\frac{1}{p}}
\geq a^{n(1-\gamma-\frac{1}{p})}\to+\infty,
\end{align*}
as $n\to+\infty$.

If $p=q=+\infty$, then it is trivial to prove that $\mathcal{N}^{\gamma}$ is unbounded on $\ell^{\infty}(\mathfrak{T})$.

Case 2: If $p=q={1}/{(1-\gamma)}$, then we take the function $g_{n,1}$ as given in Case 1 and obtain 
\[\mathcal{N}^{\gamma}g_{n,1}(x)=\mu( T(x, n-d(o,x)))^{1-\gamma},\; x\in T(o,n)\]
for every $x\in\mathfrak{T}$. It follows that 
\begin{align*}
\Vert\mathcal{N}^{\gamma}g_{n,1}\Vert_{q}^{q}
&\geq\sum_{x\in T(o,n)}\left(\mathcal{N}^{\gamma}g_{n,1}(x)\right)^{q}\\
&=\sum_{j=0}^{n}
\sum_{x\in s^{j}(o)}\mu( T(x,n-j))\\
&=\mu( T(o,n))+\sum_{j=1}^{n}(\mu( T(o,n))-\mu( T(o,j-1)))\\
&=(n+1)\mu( T(o,n))-\sum_{j=0}^{n-1}\mu( T(o,j)).
\end{align*}
It follows from Lemma \ref{BASE2} (ii) that
\begin{align*}
\sum_{j=0}^{n-1}\mu( T(o,j))
&\leq \frac{a}{a-1}\sum_{j=0}^{n-1}\mu(s^{j}(o))\\
&\leq  \frac{a}{a-1}\sum_{j=0}^{n-1}\frac{\mu(s^{n}(o))}{a^{n-j}}\\
&= \frac{a}{a-1}\mu(s^{n}(o))\sum_{k=1}^{n}\frac{1}{a^k}\\
&\leq\frac{a}{(a-1)^2}\left(1-\frac{1}{a^{n}}\right)\mu( T(o,n)),
\end{align*}
As $n\to+\infty$, we have   
\begin{align*}
\frac{\Vert\mathcal{N}^{\gamma}g_{n,1}\Vert_{q}^{q}}{\Vert g_{n,1}\Vert_{p}^{q}}
\geq n+1-\frac{a}{(a-1)^2}\left(1-\frac{1}{a^{n}}\right)\to+\infty.
\end{align*}
This completes the proof.
\end{proof}

To prove the boundedness of $\widetilde{\mathcal{N}}^\gamma$ from $\ell^{p}(\mathfrak{T})$ to $\ell^{q}(\mathfrak{T})$, we need the covering lemma of triangles on the tree. In a family of triangles in $\mathfrak{T}$, a triangle $T$ is called maximal if for any other triangle $T'$ in the family with $T' \neq T$, we have $T' \cap T \neq T$.

\begin{lem}\label{COV}
Suppose that $\mathcal{G}$ is a finite collection of maximal triangles in a locally finite tree $\mathfrak{T}$, with $3\leq a+1\leq\nu<+\infty$, and set $G:=\bigcup_{ T\in\mathcal{G}} T$. Then for every $1\leq p<+\infty$, the following estimate holds:
\[\left\Vert \sum_{ T\in\mathcal{G}}\chi_{ T}\right\Vert_{p}\leq A_{a,p}  \left\Vert\chi_{G}\right\Vert_{p},\]
where $A_{a,p}^p:=\displaystyle\frac{a^2}{a-1}\sum_{k=1}^{+\infty}\frac{k^p}{a^k}<+\infty$.
\end{lem}

The proof of this lemma follows by adapting the method of \cite[Lemma 4.1]{MS}.

\begin{rem}
For a given $r\in\mathbb{N}$ and $x\in\mathfrak{T}$, since
\[
B(x,r) = \bigcup_{j=0}^{r} T(p^{j}(x),n-j),
\]
we know that $B(x,r)$ is also a special finite collection of maximal triangles in $\mathfrak{T}$.
It then follows from Lemma \ref{COV} that
\[
\sum_{j=0}^{r}\mu(T(p^{j}(x),n-j))
\leq\frac{a^2}{a-1}\sum_{k=1}^{+\infty}\frac{k}{a^k}\mu(B(x,r))
=\Big(\frac{a}{a-1}\Big)^3\mu(B(x,r)).
\]
Proposition \ref{lem21} and Lemma \ref{COV} provide 2 different methods to prove that the volume of the ball is comparable to the sum of the volumes of a finite collection of maximal triangles.
Lemma \ref{COV} can be applied to more general collections of maximal triangles, whereas Proposition \ref{lem21} provides a sharper estimate for balls on locally finite trees.
\end{rem}

\begin{thm}\label{thm2}
Suppose that $\mathfrak{T}$ is a locally finite tree with $3\leq a+1\leq\nu<+\infty$, $0<\gamma\leq 1$, and $1<p\leq q\leq+\infty$.
Then

(i) $\widetilde{\mathcal{N}}^\gamma$ is of strong type $(p,q)$ for every $\frac{1}{q}\leq\frac{1}{p}+\gamma-1$;

(ii) $\widetilde{\mathcal{N}}^\gamma$ is unbounded from $\ell^{p}(\mathfrak{T})$ to $\ell^{q}(\mathfrak{T})$ if $\frac{1}{q}>\frac{1}{p}+\gamma-1$;

(iii) $\widetilde{\mathcal{N}}^\gamma$ is not of weak type $(1,{1}/{\gamma})$ on the homogeneous tree $\mathfrak{T}_b$, $b\geq2$.
 
\end{thm}

\begin{proof}
(i) By an argument analogous to that in \eqref{Holder}, we deduce that $\widetilde{\mathcal{N}}^\gamma$ is of strong type $({1}/{(1-\gamma)},\infty)$.
Now for any $\lambda>0$, we define the level set \begin{equation}\label{set2}
E_{\widetilde{\mathcal{N}}^\gamma f}(\lambda):=\{x\in\mathfrak{T}:\widetilde{\mathcal{N}}^\gamma f(x)>\lambda\}.
\end{equation}
If $x\in E_{\widetilde{\mathcal{N}}^\gamma f}(\lambda)$, then there exists a triangle $ T_x$ containing the vertex $x$ such that 
\[\frac{1}{\mu( T_x)^\gamma}\int_{ T_x}\vert f(y)\vert d\mu(y)>\lambda.\]
Then for every $1<p<\frac{1}{1-\gamma}$, by H\"older's inequality, we have that 
\begin{align*}
\mu( T_x)^{\gamma-1+\frac{1}{p}}
<\frac{\mu( T_x)^{1-\frac{1}{p}}}{\lambda}\int_{ T_x}\vert f(y)\vert d\mu(y)\leq\frac{\Vert f\Vert_{p}}{\lambda}.
\end{align*}
Then $E_{\widetilde{\mathcal{N}}^\gamma f}(\lambda) $ is empty if $\Vert f\Vert_{p}\leq\lambda$.
We assume $\Vert f\Vert_{p}>\lambda$. 
If $f\in\ell^p(\mathfrak{T})$, then $E_{\widetilde{\mathcal{N}}^\gamma f}(\lambda) $ can be written as a union of finite number of maximal triangles.
Denote by  $F_{\widetilde{\mathcal{N}}^\gamma f}(\lambda)$ the collection of all maximal triangles in $E_{\widetilde{\mathcal{N}}^\gamma f}(\lambda) $.
Thus, $E_{\widetilde{\mathcal{N}}^\gamma f}(\lambda) =\bigcup_{ T\in F_{\widetilde{\mathcal{N}}^\gamma f}(\lambda)} T$. 
Then for every $1<p<{1}/{(1-\gamma)}$, $0<\lambda<\Vert f\Vert_p$, and every $f$ on $\mathfrak{T}$, it follows that
\begin{align*}
\mu(E_{\widetilde{\mathcal{N}}^\gamma f}(\lambda) )^{\gamma}
&\leq\Big(\sum_{ T\in F_{\widetilde{\mathcal{N}}^\gamma f}(\lambda)}\mu( T)\Big)^\gamma\\
&\leq\sum_{ T\in F_{\widetilde{\mathcal{N}}^\gamma f}(\lambda)}\mu( T)^\gamma\\
&\leq\frac{1}{\lambda}\sum_{ T\in F_{\widetilde{\mathcal{N}}^\gamma f}(\lambda)}\int_{ T}\vert f(y)\vert d\mu(y)\\
&=\frac{1}{\lambda}\int_{\mathfrak{T}}\vert f(y)\vert\sum_{ T\in F_{\widetilde{\mathcal{N}}^\gamma f}(\lambda)}\chi_{ T}(y) d\mu(y)\\
&\leq\frac{1}{\lambda}\Vert f\Vert_{p}\left\Vert\sum_{ T\in F_{\widetilde{\mathcal{N}}^\gamma f}(\lambda)}\chi_{ T}\right\Vert_{p'}\\
&\leq\frac{A_{a,p'}}{\lambda}\Vert f\Vert_{p}\mu(E_{\widetilde{\mathcal{N}}^\gamma f}(\lambda) )^{\frac{1}{p'}},
\end{align*}
and the last inequality follows from Lemma \ref{COV}. Hence we have that 
\[\lambda\mu(E_{\widetilde{\mathcal{N}}^\gamma f}(\lambda) )^{\gamma+\frac{1}{p}-1}\leq A_{a,p'}\Vert f\Vert_{p},\]
which implies that $\widetilde{\mathcal{N}}^\gamma$ is of weak type $(p,q)$ if ${1}/{p}-{1}/{q}=1-\gamma$.
 
By the Marcinkiewicz interpolation theorem and the continuous inclusion of $\ell^{p}(\mathfrak{T})$, we then obtain that $\widetilde{\mathcal{N}}^\gamma$ is of strong type $(p,q)$ if ${1}/{p}-{1}/{q}\geq 1-\gamma$. 

(ii) Choose a triangle $ T_n$ with vertex $o$ and height $n>0$ in the locally finite tree $\mathfrak{T}$, and take the function $h_n(x)=\chi_{ T_n}(x)$, $x\in\mathfrak{T}$. Then $\Vert h_n\Vert_p=\mu(T_n)^{1/p}$.
For any vertex $x\in T_n$, we have that 
\[\widetilde{\mathcal{N}}^{\gamma}h_n(x)
\geq\frac{1}{\mu( T_n)^\gamma}\int_{ T_n}\vert h_n(y)\vert d\mu(y)
=\mu( T_n)^{1-\gamma}.\]
Thus, it follows that
\[\Big\{x\in\mathfrak{T}:\widetilde{\mathcal{N}}^{\gamma}h_n(x)
>\frac{\mu( T_n)^{1-\gamma}}{2}=:\lambda_n\Big\}\supset T_n.\]
Then for ${1}/{p}-{1}/{q}< 1-\gamma$, we have that 
\begin{align*}
\frac{\lambda_n \mu(\{x\in\mathfrak{T}:\widetilde{\mathcal{N}}^{\gamma}h_n(x)
>\lambda_n\})^{\frac{1}{q}}}{\Vert h_n\Vert_p}
\gtrsim\mu( T_n)^{\frac{1}{q}-(\frac{1}{p}+\gamma-1)}\to+\infty,
\end{align*}
as $n\to+\infty$, which proves (ii).

(iii) We now prove that $\widetilde{\mathcal{N}}^\gamma$ is not of weak type $(1,{1}/{\gamma})$ on the homogeneous tree $\mathfrak{T}_b$, $b\geq2$.
For a fixed vertex $o \in \mathfrak{T}_b$ and $\lambda > 0$, let $E_{\widetilde{\mathcal{N}}^\gamma \delta_o}(\lambda)$ be the level set given by \eqref{set2}.
Then $E_{\widetilde{\mathcal{N}}^\gamma \delta_{o}}(\lambda)$ can be written as a union of finite number of maximal triangles in which each maximal triangle with vertex $x$ and height $h$ satisfies that
\begin{align}\label{UCT}
\frac{1}{(b+1)^{\gamma}\lambda}\leq\mu( T)^{\gamma}<\frac{1}{\lambda}.
\end{align} 
The right-hand inequality follows from the definition of $E_{\widetilde{\mathcal{N}}^\gamma \delta_{o}}(\lambda)$. For the left-hand inequality, we consider a triangle $\widetilde{ T}$ with vertex $x$ and height $h+1$. Since $\widetilde{ T}\subsetneq E_{\widetilde{\mathcal{N}}^\gamma \delta_{o}}(\lambda)$ otherwise $ T$ is not the maximal triangle. It then follows that $\frac{1}{\mu(\widetilde{ T})}\leq\lambda$. Denote by $h_\lambda$ the largest integer such that the triangle $T$ with height $h_\lambda$ satisfies \eqref{UCT}. Then we have that
\[\Big(\frac{b^{h_\lambda+2}-1}{b-1}\Big)^\gamma\geq\frac{1}{\lambda},\]
which implies that 
\[b^{h_\lambda}
\geq\frac{b-1+\lambda^{\frac{1}{\gamma}}}{b^2\lambda^{\frac{1}{\gamma}}}
\geq\frac{1}{2b\lambda^{\frac{1}{\gamma}}}.\] 
Since 
\[E_{\widetilde{\mathcal{N}}^\gamma \delta_{o}}(\lambda)
\supset T(o,h_{\lambda})
\bigcup\left(\bigcup_{k=1}^{h_\lambda}
( T(p^k(o),h_{\lambda})\backslash T(p^{k-1}(o),h_{\lambda}))\right)\]
and
\[\mu( T(p^k(o),h_{\lambda})\backslash T(p^{k-1}(o),h_{\lambda}))
=\frac{b^{h_{\lambda}+1}-1}{b-1}-\frac{b^{h_{\lambda}}-1}{b-1}=b^{h_\lambda},\]
for $1\leq p\leq q\leq {1}/{\gamma}$, it follows that 
$\mu(E_{\widetilde{\mathcal{N}}^\gamma \delta_{o}}(\lambda))
\geq h_\lambda b^{h_\lambda}
\geq\frac{\log_b \frac{1}{2b\lambda^{1/\gamma}}}{2b\lambda^{1/\gamma}}$. Let $\eta:=\lambda^{\frac{1}{\gamma}}$ and $t:=\log_b \frac{1}{2b\eta}$. Then
\begin{align*}
\frac{\lambda\mu(E_{\widetilde{\mathcal{N}}^\gamma \delta_{o}}(\lambda))^{\frac{1}{q}}}{\Vert\delta_o\Vert_p}
\gtrsim\frac{t^{1/q}}{b^{t(\gamma-1/q)}}\to+\infty,\quad t\to+\infty.
\end{align*}
This completes the proof.
\end{proof}

\begin{figure}[htbp]
  \centering

  \begin{subfigure}[c]{0.45\textwidth}
    \centering
    \begin{tikzpicture}[scale=3.8]
      \def\gammaVal{0.7}

      \draw[->] ({1-\gammaVal},0) -- (1.1,0) node[right] {$1/p$};
      \draw[->] (0,0) -- (0,1.1) node[above] {$1/q$};

      \coordinate (O) at (0,0);
      \coordinate (A) at ({1-\gammaVal},{1-\gammaVal});
      \coordinate (B) at (1,1);
      \coordinate (C) at (1,0);

      \fill[pattern= vertical lines, pattern color=black!30]
        (O) -- (A) -- ({1-\gammaVal},0) -- cycle; 

      \fill[pattern=horizontal lines, pattern color=black!30]
        (A) -- (B) -- (C) -- ({1-\gammaVal},0) -- cycle; 
      
      \draw[black!100,thick] (A) -- (B) -- (C) -- ({1-\gammaVal},0);
      \draw[loosely dashed, black!100, thick]  (O)--({1-\gammaVal},0) -- cycle;
      \draw[loosely dashed, black!100, thick] (O) -- (A) -- cycle ;
      
      \draw[dotted, black,very thick] (A) -- ({1-\gammaVal},0);

      \draw[black, line width=0.2pt] (O) circle (0.7pt);  
      \fill[white] (O) circle (0.7pt);                 
      \fill[black] (O) circle (0.2pt);                 
      
      \fill[black!100] ({1-\gammaVal},0) circle (0.6pt) node[below, text=black] {$1-\gamma$};
      \fill[black!100] (C) circle (0.6pt) node[below, text=black] {$1$};
      
      \draw[black, line width=0.2pt] (A) circle (0.7pt);  
      \fill[white] (A) circle (0.7pt);                 
      \fill[black] (A) circle (0.2pt);                 
  
      \fill[black!100] (B) circle (0.6pt);

      \draw[dashed] (0,{1-\gammaVal}) -- ({1-\gammaVal},{1-\gammaVal});
      \draw[dashed] (0,1) -- (1,1);

      \node[below] at (0,0) {$0$};
      \node[left] at (0,{1-\gammaVal}) {$1-\gamma$};
      \node[left] at (0,1) {$1$};
    \end{tikzpicture}
    \caption{$\mathcal{N}^\gamma$}
  \end{subfigure}
  \hfill
  \begin{subfigure}[c]{0.45\textwidth}
    \centering
    \begin{tikzpicture}[scale=3.8]
      \def\gammaVal{0.7}

      \draw[->] ({1-\gammaVal},0) -- (1.1,0) node[right] {$1/p$};
      \draw[->] (0,0) -- (0,1.1) node[above] {$1/q$};

      \coordinate (O) at (0,0);
      \coordinate (A) at ({1-\gammaVal},0);
      \coordinate (B) at (1,1);
      \coordinate (C) at (1,0);
      \coordinate (D) at (1,{\gammaVal});
      
      \fill[pattern=vertical lines, pattern color=black!30]
        (O) -- (A) -- (D) -- (B) -- cycle;

      \fill[pattern=horizontal lines,pattern color=black!30]
        (A) -- (C) -- (D) -- cycle;
      
      \draw[black!100,thick] (A) -- (D);
      \draw[black!100,thick] (A) -- (C);
      \draw[dotted, black!100, very thick] (D) -- (C);
      \draw[loosely dashed, black!100,thick]  (O) -- (B);
      \draw[loosely dashed, black!100,thick]  (D) -- (B);
      \draw[loosely dashed, black!100,thick]  (O) -- (A);

      \draw[black, line width=0.2pt] (O) circle (0.7pt);  
      \fill[white] (O) circle (0.7pt);                 
      \fill[black] (O) circle (0.2pt);                 
      
      \draw[black, line width=0.3pt] (C) circle (0.7pt);  
      \fill[white] (C) circle (0.7pt);                 
      \draw[black, line width=0.3pt] (C) circle (0.4pt)  
      node[below, text=black] {$1$};
      
      \draw[black, line width=0.2pt] (B) circle (0.7pt);  
      \fill[white] (B) circle (0.7pt);                 
      \fill[black] (B) circle (0.2pt);                 
      
      \fill[black!100]  (A) circle (0.6pt) node[below, text=black] {$1-\gamma$};
      
      \draw[black, line width=0.2pt] (D) circle (0.7pt);  
      \fill[white] (D) circle (0.7pt)node[above right, text=black]{};                    
               
      \draw[dashed] (0,1) -- (1,1);
      \draw[dashed] (0,{\gammaVal}) -- (1,{\gammaVal});

      \node[below] at (0,0) {$0$};
      \node[left] at (0,1) {$1$};
      \node[left] at (0,{\gammaVal}) {$\gamma$};
      
    \end{tikzpicture}
    \caption{$\widetilde{\mathcal{N}}^\gamma$}
  \end{subfigure}

\caption{The $\ell^p-\ell^q$ bounded regions of the maximal operators $\mathcal{N}^\gamma$ and $\widetilde{\mathcal{N}}^\gamma$ $(0<\gamma<1)$ on locally finite trees in the $(1/p,1/q)$ plane.}
\end{figure}
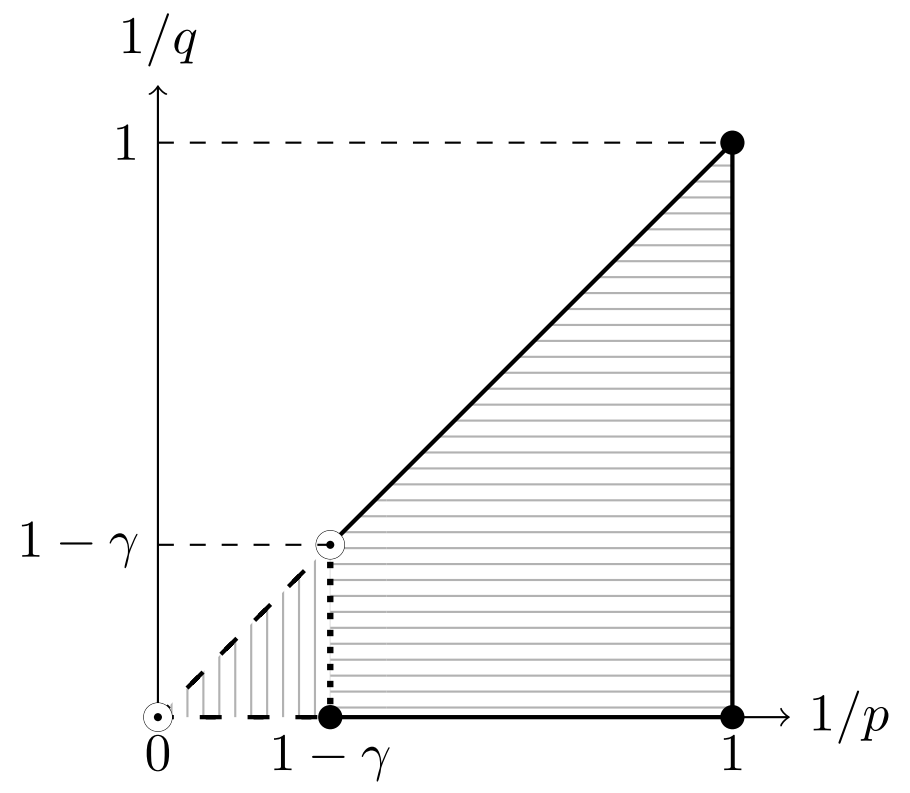
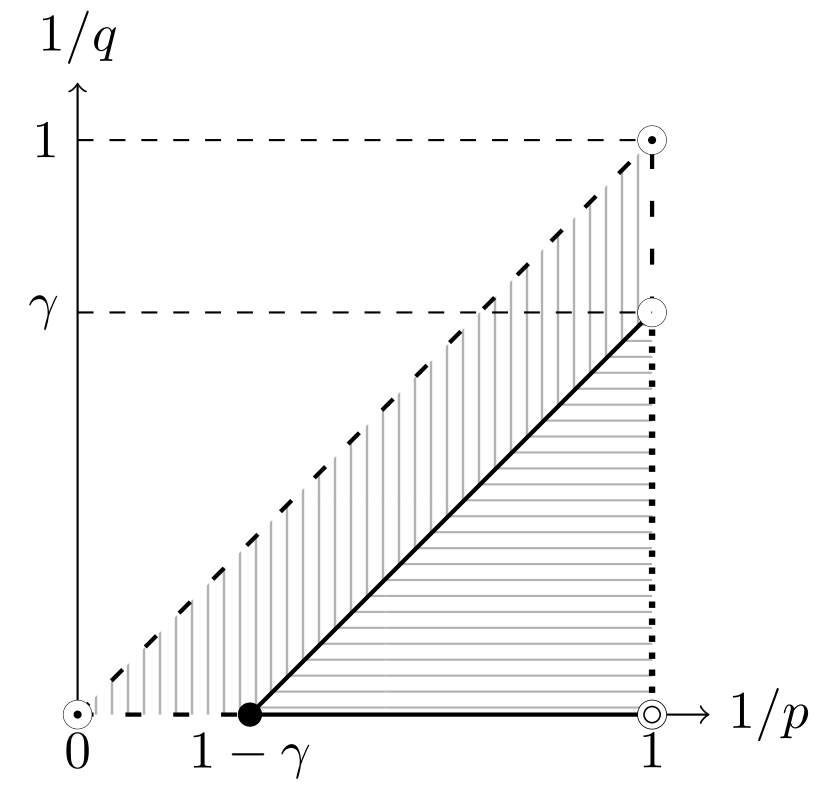

The figure above summarizes all the results that we obtained on the strong types and weak types boundedness of $\mathcal{N}^\gamma$ and $\widetilde{\mathcal{N}}^\gamma$. 
The region filled with horizontal lines, including solid lines and points, represents parameter pairs where the operators are of strong type $(p,q)$.
The region filled with vertical lines, including dashed lines and hollow points with a dot inside, indicates where the operators fail to be strong type $(p,q)$, with the hollow point specifically denoting failure of weak type $(p,q)$.
The dotted lines and hollow point with circle inside is the parameter range where boundedness remains undetermined.

\begin{rem}
Compared to \cite[Theorem 3.4 and Fig. 1]{LS}, we see that $\mathcal{N}^\gamma$ has a larger range of $p$ and $q$ while $\widetilde{\mathcal{N}}^\gamma$ has a small range of $p$ and $q$ than $\mathcal{M}^\gamma$ to be of strong type $(p,q)$ on the homogeneous tree $\mathfrak{T}_b$. 
\end{rem}

\section{Final Remarks: Unboundedness of $\widetilde{\mathcal{M}}^\gamma$ in Contrast to $\widetilde{\mathcal{N}}^\gamma$ on $\mathfrak{F}_{a,b}$}\label{s5}

Proposition \ref{p3} (ii) establishes that ${\mathcal{N}}^\gamma \lesssim \widetilde{\mathcal{M}}^\gamma$ on $\Upsilon_{a,b}$ for $\gamma \geq \frac{1}{2\log_b a - 1}$ and $a<b<a^2$, leaving open the behavior for smaller $\gamma$. Our main contribution in this final section is to resolve this gap for a key example. We prove that on the specific tree $\mathfrak{F}_{a,b} \in \Upsilon_{a,b}$, and for all $0<\gamma\leq 1$, the ball-type maximal operator $\widetilde{\mathcal{M}}^\gamma$ is indeed substantially less bounded than $\widetilde{\mathcal{N}}^\gamma$, exhibiting unboundedness on a strictly larger set of $(p,q)$ pairs. 

\begin{thm}\label{thm51}
Suppose that $0 < \gamma \leq 1$, $1 \leq p \leq q < +\infty$, and $b > a \geq 2$. Then $\widetilde{\mathcal{M}}^\gamma$ is unbounded from $L^p(\mathfrak{F}_{a,b})$ to $L^q(\mathfrak{F}_{a,b})$ if $\frac{1}{q}>\frac{1}{2-\log_b a}\Big(\frac{1}{p}+\gamma-1\Big)$. 
\end{thm}
Before the proof of the theorem, we recall that a tree $\mathfrak{T}$ is said to satisfy the Cheeger isoperimetric property if there exists a constant $C>0$ such that for every subset $E \subseteq \mathfrak{T}$,
\[
\mu(\partial E)\geq C\mu( E),
\]
where $\partial E$ denotes the boundary of $E$,
that is, the collection of vertices in $E$ having at least one adjacent vertex in $\mathfrak{T}\setminus E$. 
The supremum of all such constants $C$ is called the Cheeger isoperimetric constant of $\mathfrak{T}$, denoted by $C_{\mathfrak{T}}$. 

\begin{prop}[\cite{RT}]\label{p51}
Let $2\leq a\leq b<+\infty$ and $\mathfrak{T}$ be a tree in the class $\Upsilon_{a,b}$. Then $\mathfrak{T}$ satisfies the Cheeger isoperimetric property.
\end{prop}

For any $x\in\mathfrak{F}_{a,b}$ and $r\in\mathbb{N}$, denote by $S(x,r):=\{y\in\mathfrak{T}:d(x,y)=r\}$. Since $\mathfrak{F}_{a,b}$ is in the class $\Upsilon_{a,b}$, it  follows from Proposition \ref{p51} that 
\begin{align}\label{SPH}
\mu(S(x,r))
\leq \mu(B(x,r))
\leq C_{\mathfrak{F}_{a,b}}^{-1}\mu(S(x,r)).
\end{align}

Since the measure $\mu(S(x,r))$ is more tractable than $\mu(B(x,r))$ on the tree $\mathfrak{F}_{a,b}$, 
the equivalence \eqref{SPH} provides a more direct route to estimate the maximal operator $\widetilde{\mathcal{M}}^\gamma$.

\begin{proof}[Proof of Theorem \ref{thm51}.]
Let $x_0$ be the reference vertex of $\mathfrak{F}_{a,b}$. For large integers $n$, denote by $\tau:=\log_a b$ and the set
\[
E_n
:=\Big\{y\in\mathfrak{F}_{a,b}\backslash\bigcup_{k=0}^{+\infty}s^{k}(x_0): d(y,x_0)=n\Big\}
\] 
and
\[
F_n
:=\Big\{y\in\bigcup_{k=0}^{+\infty}s^{k}(x_0): d(y,x_0)=\lfloor(2\tau-1)n\rfloor\Big\}.
\]
For $y\in F_n$, there is a unique vertex $z$ on the geodesic joining $y$ and $x_0$ at distance $\lfloor\tau n\rfloor$ from $y$. Since 
\[
d(x_0,z)=\lfloor(2\tau-1)n\rfloor-\lfloor\tau n\rfloor
<(2\tau-1)n-\tau n+1=(\tau-1)n+1,
\]
it follows that for any $x\in E_n$,
\[
d(z,x)=d(z,x_0)+d(x_0,x)
=\lfloor(2\tau-1)n\rfloor-\lfloor\tau n\rfloor+n
<\tau n+1.
\]
Choose $r_n=\lfloor\tau n\rfloor+1$. Then the ball $B(z,r_n)$ contains both $y$ and whole $E_n$. Hence, given $x\in F_n$, and by \eqref{SPH} there exists a vertex $z\in\mathfrak{F}_{a,b}$ such that 
\[
\widetilde{\mathcal{M}}^\gamma\chi_{E_n}(x)
\geq\frac{1}{\mu(B(z,r_n))^\gamma}\int_{B(z,r_n)}\chi_{E_n}(y)d\mu.
\] 
It is proved in \cite[Theorem 4.1 (i)]{LMSV} that 
\begin{align*}
\mu(B(z,r_n))\leq C_{\mathfrak{F}_{a,b}}^{-1}\mu(S(z,r_n))
&=a^{r_n}+\sum_{k=1}^{d(z,x_0)}(a-1)a^{r_n-k-1}+b^{r_n-d(z,x_0)-1}
\\
&=(a+1)a^{r_n-1}+b^{r_n-d(z,x_0)-1}-a^{r_n-d(z,x_0)-1}.
\end{align*} 

Since
\begin{align*}
r_n-d(z,x_0)-1
&=\lfloor\tau n\rfloor+1-(\lfloor(2\tau-1)n\rfloor-\lfloor\tau n\rfloor)\\
&<2\tau n-((2\tau-1)n-1)+1\\
&=n+2{\Large }
\end{align*}
and 
\begin{align*}
r_n-d(z,x_0)-1
>2(\tau n-1)+1-(2\tau-1)n=n-1,
\end{align*}
it follows that
\begin{align*}
\mu(S(z,r_n))<(a+1)a^{\tau n}+b^{n-2}-a^{n-1}<(a+1+b^2)b^n.
\end{align*}

There is another method to estimate the measure of $B(z,r_n)$. In fact, by Lemma \ref{BASE2} (ii) and \eqref{MUBUT}, we have that 
\begin{align*}
\mu(B(z,r_n))
&\leq\sum_{j=0}^{r_n}\mu(T(p^j(z),r_n-j))\\
&\leq\frac{a}{a-1}\sum_{j=0}^{r_n}\mu(\beta(T(p^j(z),r_n-j)))\\
&=\frac{a}{a-1}\left(\sum_{j=0}^{d(z,x_0)}a^{r_n-j}
+\sum_{j=d(z,x_0)+1}^{r_n}\Big((b-1)b^{r_n-j-1}+a^{r_n-j-1}\Big)\right)\\
&=\frac{a}{a-1} \left( \frac{a^{r_n+1}-a^{-1}}{a-1}-a^{r_n-d(z,x_0)-1} + b^{r_n-d(z,x_0)-1}-\frac{1}{b} \right)\\
&<\frac{a}{a-1}\Big(\frac{a^{\tau n+2}}{a-1}-a^{n-1}+b^{n+2}-\frac{1}{a(a-1)}-\frac{1}{b}\Big)\\
&=\frac{a}{a-1}\Big(\frac{a^2}{a-1}+b^2\Big)b^n-\frac{a^n}{a-1}-\frac{a^{2} - a + b}{b (a-1)^{2}}
\end{align*}  
and
\[
\mu(F_n)=a^{\lfloor(2\tau-1)n\rfloor}
\geq a^{(2\tau-1)n-1}=\frac{1}{a}\Big(\frac{b^2}{a}\Big)^n.
\]
Thus, we can conclude that for $1\leq p\leq q<+\infty$,
\begin{align*}
\frac{\Vert\widetilde{\mathcal{M}}^\gamma\chi_{E_n}\Vert_q^q}{\Vert \chi_{E_n}\Vert_p^q}
&\geq
\frac{\sum_{x\in F_n}\widetilde{\mathcal{M}}^\gamma\chi_{E_n}(x)^q}{\mu(E_n)^{q/p}}\\
&\geq\frac{\mu(F_n)\mu(E_n)^q}{\mu(E_n)^{q/p}\mu(B(z,r_n))^{q\gamma}}\\
&\gtrsim\frac{(b^2/a)^n b^{nq(1-\gamma)}}{b^{nq/p}}\\
&=b^{(2+q(1-\gamma)-1/\tau-q/p)n}.
\end{align*}
Therefore if $2+q(1-\gamma)-1/\tau-q/p>0$, that is, $\frac{1}{q}>\frac{1}{2-\log_b a}\Big(\frac{1}{p}+\gamma-1\Big)$, then 
\[
\frac{\Vert\widetilde{\mathcal{M}}^\gamma\chi_{E_n}\Vert_q^q}{\Vert \chi_{E_n}\Vert_p^q}
\to+\infty,
\]
as $n\to+\infty$. This completes the proof.
\end{proof}

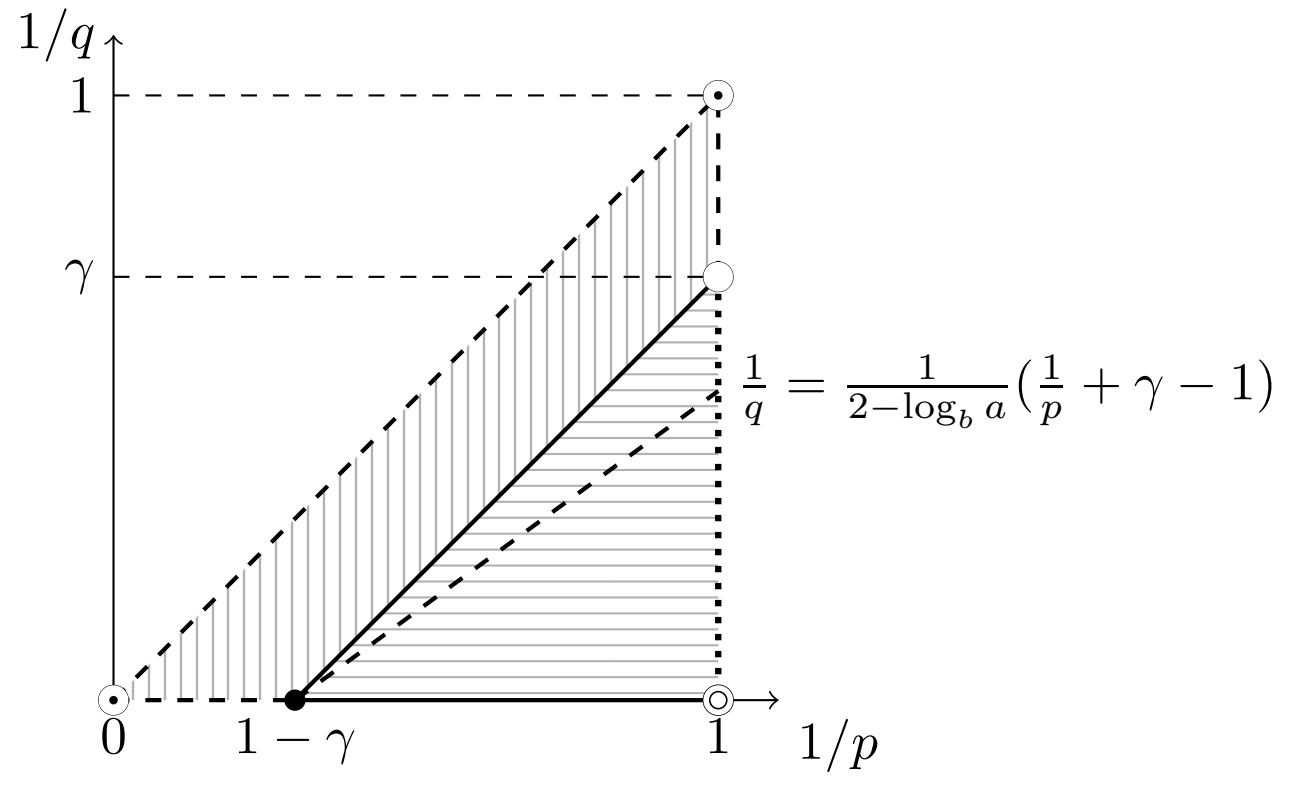
\begin{figure}[htbp]
    \centering
    \begin{tikzpicture}[scale=4]
      \def\gammaVal{0.7}
      
      \pgfmathsetmacro{\logValue}{ln(2)/ln(3)}  
      \pgfmathsetmacro{\slopeValue}{1/(2-\logValue)}  
      
      \draw[->] ({1-\gammaVal},0) -- (1.1,0) node[below right] {$1/p$};
      \draw[->] (0,0) -- (0,1.1) node[left] {$1/q$};

      \coordinate (O) at (0,0);
      \coordinate (A) at ({1-\gammaVal},0);
      \coordinate (B) at (1,1);
      \coordinate (C) at (1,0);
      \coordinate (D) at (1,{\gammaVal});
      
      \fill[pattern=vertical lines, pattern color=black!30] 
        (O) -- (A) -- (D) -- (B) -- cycle;

      \fill[pattern=horizontal lines, pattern color=black!30]
        (A) -- (C) -- (D) -- cycle;
      
      \draw[black,thick] (A) -- (D);
      \draw[black,thick] (A) -- (C);
      
      \draw[black,very thick, dotted] (D) -- (C);
      
      \draw[black,thick, dashed] (O) -- (B);
      \draw[black,thick, dashed] (D) -- (B);
      \draw[black,thick, dashed] (O) -- (A);
      
      \pgfmathsetmacro{\yAtOne}{\gammaVal/(2-\logValue)}  
      
      \coordinate (E) at ({1-\gammaVal},0);  
      \coordinate (F) at (1,{\yAtOne});      
      
      \draw[black, thick, dashed] (E) -- (F);
      
      \node[right, text=black] at (F) {$\frac{1}{q}=\frac{1}{2-\log_b a}(\frac{1}{p}+\gamma-1)$};

      \draw[black, line width=0.2pt] (O) circle (0.7pt);  
      \fill[white] (O) circle (0.7pt);                 
      \fill[black] (O) circle (0.2pt);                 
      
      \draw[black, line width=0.2pt] (B) circle (0.7pt);  
      \fill[white] (B) circle (0.7pt);                 
      \fill[black] (B) circle (0.2pt);                 
      
      \draw[black, line width=0.3pt] (C) circle (0.7pt);  
      \fill[white] (C) circle (0.7pt);                 
      \draw[black, line width=0.3pt] (C) circle (0.4pt);  
      \node[below, text=black] at (C) {$1$};
      
      \fill[black] (A) circle (0.5pt) node[below, text=black] {$1-\gamma$};
      
      \draw[black, line width=0.2pt] (D) circle (0.7pt);  
      \fill[white] (D) circle (0.7pt);                    
      \node[above right, text=black] at (D) {};

      \draw[dashed] (0,1) -- (1,1);
      \draw[dashed] (0,{\gammaVal}) -- (1,{\gammaVal});

      \node[below] at (0,0) {$0$};
      \node[left] at (0,1) {$1$};
      \node[left] at (0,{\gammaVal}) {$\gamma$};   
    \end{tikzpicture}
\caption{The function $\frac{1}{q}=\frac{1}{2-\log_b a}(\frac{1}{p}+\gamma-1)$, $(2\leq a<b)$ and the $(\ell^p,\ell^q)$ bounded regions of $\widetilde{\mathcal{N}}^\gamma$ $(0<\gamma< 1)$ on $\mathfrak{F}_{a,b}$ in the $(1/p,1/q)$ plane.}
\end{figure}

Hence, for the range $0<\gamma\leq 1$, the question of whether $\widetilde{\mathcal{M}}^\gamma$ pointwise dominates $\widetilde{\mathcal{N}}^\gamma$ on every locally finite tree is still open. However, on the specific tree $\mathfrak{F}_{a,b}$ ($2<a<b$), we have conclusively shown that the operator $\widetilde{\mathcal{M}}^\gamma$ exhibits a strictly wider region of unboundedness, as defined by $\frac{1}{q} > \frac{1}{2-\log_b a}(\frac{1}{p}+\gamma-1)$, compared to $\widetilde{\mathcal{N}}^\gamma$, whose unboundedness occurs when $\frac{1}{q} > \frac{1}{p}+\gamma-1$.

\end{document}